\documentclass[11pt]{article}%
\usepackage{amsmath}%
\usepackage{amsfonts}%
\usepackage{amssymb}%
\usepackage{graphicx}
\providecommand{\U}[1]{\protect\rule{.1in}{.1in}}
\newcommand{\sfrac}[2]{\frac{\displaystyle #1}{\displaystyle #2}}
\newcommand{\ignore}[1]{}
\newcounter{subscript}

\def\Complex{\rm{C\!\!\!I}}

\begin{document}
\title{Estimating spectral density functions for Sturm-Liouville problems with two
singular endpoints}
\date{  }
\author{CHARLES FULTON\\Department of Mathematical Sciences \\Florida Institute of Technology \\Melbourne, Florida 32901-6975 \\ \\
DAVID PEARSON \\Department of Mathematics \\University of Hull \\Cottingham Road \\Hull, HU6 7RX England \\ \\
STEVEN PRUESS \\1133 N Desert Deer Pass \\Green Valley, Arizona 85614-5530
}
\maketitle


\abstract{In this paper we consider the Sturm-Liouville equation $-y'' + qy = \lambda y$ on the
half line $(0,\infty)$ under the assumptions that $x=0$ is a regular singular point and nonoscillatory
for all real $\lambda$, and that either (i) $q$ is $L_1$ near $x=\infty$, or (ii)  $q'$ is $L_1$ near $\infty$
with $q(x) \rightarrow 0$ as $x \rightarrow \infty$, so that there is absolutely continuous spectrum
in $(0,\infty)$.  Characterizations of the spectral density function for this doubly singular problem,
similar to those obtained in \cite{FPP2} and \cite{FPP4}
(when the left endpoint is regular) are established;
corresponding approximants from the two algorithms in \cite{FPP2} and \cite{FPP4}
are then utilized, along
with Frobenius recurrence relations and piecewise trigonometric - hyperbolic splines, to generate numerical
approximations to the spectral density function associated with the doubly singular problem on $(0,\infty)$.
In the case of the radial part of the separated hydrogen atom problem, the new algorithms are capable of
achieving near machine precision accuracy over the range of $\lambda$ from 0.1 to 10000, accuracies which could not
be achieved using the SLEDGE software package.}

\vskip 0.4in {\footnotesize This research partially supported by National
Science Foundation Grant DMS-0109022 to Florida Institute of Technology, and
by Engineering and Physical Sciences Research Council Grant S63403/01 to University of Hull.}

\pagebreak

\section{Introduction}

\setcounter{equation}{0}

In this paper we consider the Sturm-Liouville equation,
\begin{equation}
-y^{\prime\prime}+q(x)y=\lambda y,\label{1.1}%
\end{equation}
on $(0,\infty)$, with two singular endpoints on the half line $(0,\infty)$
under the same assumptions as in  \cite{MN,FL} under which $x=0$ is a
regular singular point of type {\bf LC/N} or {\bf LP/N} (limit circle or limit point
and nonoscillatory at $x=0$ for all real $\lambda$) and $x=\infty$ is of type
{\bf LP/O-N} with cutoff $\Lambda=0$ (limit point and nonoscillatory at $x=\infty$
for $\lambda\in(-\infty,0)$ and oscillatory for $\lambda\in(0,\infty)$);
see \cite[p. 114]{TOMS98} for these definitions of 
endpoint classifications developed in connection with the \textbf{SLEDGE}
software package. Under these assumptions, the spectrum is simple and the eigenfunction expansion
associated with equation (\ref{1.1}) has the general form

\bigskip%
\begin{align}
f(x)  & =\int_{-\infty}^{\infty}T(\lambda)\cdot\phi(x,\lambda)d\rho
(\lambda)\nonumber \\
& =\sum_{\lambda_{n}\leq0}\left[  \frac{\int_{0}^{\infty}f\left(  t\right)
\phi\left(  t,\lambda_{n}\right)  dt}{\left\Vert \phi\left(  t,\lambda
_{n}\right)  \right\Vert ^{2}}\right]  \phi\left(  x,\lambda_{n}\right)
+\int_{0}^{\infty}T\left(  \lambda\right)  \cdot\phi\left(  x,\lambda\right)
d\rho\left(  \lambda\right) \label{1.2}
\end{align}
where%

\begin{equation}
T(\lambda):=\lim_{b\rightarrow\infty}\int_{0}^{b}f(x)\phi(x,\lambda
)dx,\label{1.3}%
\end{equation}
and the solution $\phi(\cdot,\lambda)$ is a suitably normalized Frobenius
solution near the regular singular point $x=0$.

If we make, in addition to the above assumptions, the more stringent assumptions
posed in \cite{FPP2} (see {\bf Assumption 3} below), then there is absolutely continuous (a.c.)
spectrum  in $(0,\infty)$, and the spectral function $\rho(\lambda)$ in
(\ref{1.2}) is absolutely continuous on all closed intervals in $(0,\infty)$.

The purpose of this paper is \textbf{(i)} \ to extend the analysis in
\cite{FPP2,FPP4} \ under suitable assumptions on \ q(x) \ to \ show that the
spectral density function associated with (\ref{1.2})\ over the \ a.c. range $(0,\infty)$, \ that is, 
\vspace{0.1in}

\begin{equation}
f(\lambda)=\rho^{\prime}(\lambda),\ \ \ \ \ \ \ \rho(\lambda)=\int
_{0}^{\lambda}f(\mu)\ \ d\mu\bigskip\label{1.4}%
\end{equation}
\vspace{0.1in}
\noindent
can be represented for all \ $\lambda\in(0,\infty)$ \ as (see {\bf Theorem 3} below) 
\vspace{0.1in}

\begin{equation}
f(\lambda)=\frac{\displaystyle1}{\displaystyle\pi\lbrack P(x,\lambda
)\phi(x,\lambda)^{2}+Q(x,\lambda)\phi(x,\lambda)\phi^{\prime}(x,\lambda
)+R(x,\lambda)\phi^{\prime}(x,\lambda)^{2}]}.\label{1.5}%
\end{equation}
\vspace{0.1in}
\noindent
where \ $(P(\cdot,\lambda),Q(\cdot,\lambda),R(\cdot,\lambda))^{T}$ \ is the
unique solution of the initial \ value problem at \ $x=\infty:$

\begin{eqnarray}
\frac{dU}{dx}=\frac{d}{dx}\left[
\begin{array}
[c]{c}%
P\\
Q\\
R
\end{array}
\right] &=&\left[
\begin{array}
[c]{ccc}%
0 & \lambda-q & 0\\
-2 & 0 & 2(\lambda-q)\\
0 & -1 & 0
\end{array}
\right]  \cdot\left[
\begin{array}
[c]{c}%
P\\
Q\\
R
\end{array}
\right]  .\label{1.6}
\\
\lim_{x\rightarrow\infty}\left(
\begin{array}
[c]{l}
P(x,\lambda)\\
Q(x,\lambda)\\
R(x,\lambda)
\end{array}
\right) &=&\left(
\begin{array}
[c]{l}%
\sqrt{\lambda}\\
0\\
\frac{1}{\sqrt{\lambda}}%
\end{array}
\right)  ,\text{ \ \ \ }\lambda\in(0,\infty),\label{1.7}%
\end{eqnarray}
\vspace{0.1in}
and \ \textbf{(ii)} \ to extend the numerical algorithms from
\ \cite{FPP2,FPP4} \ for the computation of the spectral density function $f(\lambda)$.
% \ Here \ the right hand side of \ref{1.5} is independent of
%\ $x\in(0,\infty).$ \ In contrast to \cite{FPP2,FPP4} \ the solution
%\ $\phi(\cdot,\lambda),$ \ in terms of which the eigenfunction expansion is
%expressed, is not normalized by standard initial conditions at a regular
%point, but instead by normalization of a Frobenius \ solution of \ \ref{1.1}
%at the \textbf{R.S.P \ }$x=0.$
We illustrate the new numerical algorithms on
several examples, including \ the \ radial part of the separated hydrogen
atom. For the first objective {\bf (i)} we make use of the fact that $f(\lambda)$ is
characterized as the boundary value of a suitable Titchmarsh-Weyl m-function by
the Titchmarsh-Kodaira formula, which was recently established for such doubly 
singular problems in \cite{GZ,FL} (see {\bf Theorem 2} below). For the second
objective {\bf (ii)} we make use of exact Frobenius power series to estimate the
solution $\phi(x,\lambda)$ and its derivative near $x=0$, and then apply  initial
conditions at a suitable point $x_{0}(\lambda) > 0$ using values of
$\phi(x_{0}(\lambda),\lambda)$ and $\phi^{\prime}(x_{0}(\lambda),\lambda)$ which can generally be computed
to machine precision, so that 
numerical algorithms from \cite{FPP2,FPP4} with the left endpoint regular can be adapted
to approximate the right hand side of (\ref{1.5}); this is done by shooting with
piecewise trigonometric / hyperbolic splines to compute the solution $\phi$ of (\ref{1.1})
and the solution $(P,Q,R)^T$ of (\ref{1.6}) at a suitable `matching' point
$x \in (x_{0}(\lambda),\infty)$.

The organization of topics needed in this paper to accomplish the above two objectives
is as follows:  In section 2 we give the main assumptions near $x=0$ from \cite{MN,FL},
and the general forms of two linearly independent Frobenius power series solutions in 
all the cases we consider in this paper. In section 3 we list (without proof) the elementary
results which relate and interconnect solutions of the Sturm-Liouville equation (\ref{1.1})
with solutions of Appell's first order system (\ref{1.6}); none of these elementary results
require any special assumptions on the potential $q$. In section 4 we add the main assumptions
from \cite{FPP2} under which the initial value problem (\ref{1.6})-(\ref{1.7}) has a unique
solution in $(0,\infty)$, and reformulate in terms of solutions of Appell's system (\ref{1.6})
results obtained by D.B. Pearson and his student Al-Naggar in \cite{AlN1,AlN2}. This yields
the spectral density function characterization (\ref{1.5}) in the relatively simple case of a regular 
left endpoint. In contrast to \cite{AlN1,AlN2} we do not focus the analysis on the 
third order ordinary differential equation (see (\ref{4.13}) below) which is
satisfied by the third component, $R(x,\lambda)$, of (\ref{1.6}), but make use instead of
many of the elegant formulas from section 2. In section 5 we generalize the methods of
\cite{AlN1,AlN2} to the doubly singular problem on $(0,\infty)$, when the m-function
is defined relative to the suitably normalized Frobenius fundamental system as in \cite{MN,FL}.
This yields {\bf Theorem 3} below (a new result) in which the new Titchmarsh-Kodaira formula (see {\bf Theorem 2} below) for
the (single) spectral density function associated with the doubly singular problem gets
converted to the form (\ref{1.5}). In section 6 we list four test examples of equations on $(0,\infty)$
from \cite{MN,FL} for which explicit closed form formulas for the spectral density function were
obtained. In section 7 we describe how the two different types of numerical algorithms from
\cite{FPP2} and \cite{FPP4} (for cases involving a regular left endpoint) can be adapted, using
appropriate heuristics, to yield new algorithms for computing the spectral density function
when both endpoints are singular, which is done by utilizing the characterization (\ref{1.5})
in an appropriate way. In section 8 we give numerical output showing that the new algorithms
for doubly singular problems can achieve very high accuracy on the four test examples in section 6; we also
give numerical ouput demonstrating convergence of our numerical approximations for a potential on $(0,\infty)$
from quantum chemistry where the potential has an infinite series representation satisfying all
our assumptions. In section 9 we make use of our new code, AutoB, for the spectral density computation
in order to generate, by quadratures, approximations to the spectral functions for the four examples
in section 6, and give comparisons on timing and accuracy with the corresponding {\bf SLEDGE} runs.
Our main conclusion is that the new algorithms for doubly singular problems are very much superior
to the older algorithms for doubly singular problems which were implemented in the {\bf SLEDGE} software package. 

\noindent {\bf  Remark.}      In our
previous papers  \cite{FPP2,FPP4}  the system (\ref{1.6})  was
referred to as the \textquotedblleft PQR equations"  (our notation);
 however,  the analysis leading to them (particularly \ the motivating
 property  (\ref{3.15}) )  was discovered by M. Appell  \cite{APPELL}
in  1880. Accordingly, we will henceforth refer to this first order system
 as  the  Appell  equations.

\section{Suitably Normalized Frobenius Solutions}

\setcounter{equation}{0}

In this section we repeat the basic definitions and some of the elementary
properties of the suitably normalized Frobenius solutions which were
introduced in Fulton \cite{MN} and Fulton and Langer \cite{FL}. It will be
noted that, in most of the common cases (Bessel quations, Confluent
Hypergeometric equations, Whittaker equations) rather standard normalizations
of well known special functions have to be abandoned in order to achieve the
desired analytic properties of the Frobenius solutions (particularly, entire behaviour
in $\lambda$) needed for the fundmental definition of a single Titchmarsh-Weyl
m-function, the corresponding scalar spectral function, and for determination
of eigenfunction expansions of the problems considered in this paper (all of
which have simple spectrum).

We consider in this paper the Sturm-Liouville equation (\ref{1.1}) on the half
line $(0,\infty)$ under the following assumptions (see \cite{MN,FL}):

\begin{center}
{\bf Assumption 1:} $\mathbf{Near\ x=0:}$%
\end{center}
\underline{\textbf{Case I:}} \ For all \ $x\in(0,\infty),$%

\begin{equation}
q(x)=\frac{q_{0}}{x^{2}}+\frac{q_{1}}{x^{{}}}+\sum_{n=0}^{\infty}q_{n+2}%
x^{n},\text{ \ }q_{n}\text{ \ real \ for all \ }n,\label{2.1}%
\end{equation}
where the series is convergent in $(0,\infty)$, and where
\begin{equation}
-\frac{1}{4}\leq q_{0}<\infty.\text{ \ and \ }q_{0},q_{1}\text{ \ \ not both
\ zero.}\label{2.2}%
\end{equation}
or

\noindent\underline{\textbf{Case II:}} There exists $a>0$ such that $q(x)$ is
given by (\ref{2.1}) for $x\in(0,a]$ where the series is convergent in $(0,a]$
and where (\ref{2.2}) holds, and in the interval $[a,\infty)$ we have $q\in
L_{loc}^{1}[a,\infty)$,

\bigskip

and%

\begin{center}
{\bf Assumption 2:} $\mathbf{Near\ x=\infty:}$%
\end{center}

\begin{equation}
\lim_{x\rightarrow\infty}q(x)=0,\label{2.3}%
\end{equation}
The assumptions near $x=0$ ensure that the indicial roots near the regular
singular point $x=0$ are both real; it follows that the endpoint $x=0$ is
either \textbf{LC/N or LP/N}, and the assumption near $x=\infty$ ensure that
the endpoint $x=\infty$ is \textbf{LP/O-N} with cutoff $\Lambda = 0$ in the terminology of
\cite{TOMS98}.  Under the above assumptions it was \ proved in \cite[Theorems
4.2,4.3,5.3,5.4]{MN} \ and \ \cite[Theorem 4.5]{FL} \ that the
eigenfunction expansion associated with (\ref{1.1}) \ (in both the\textbf{ LP
\ }$and$\textbf{ \ LC} cases at $x=0$ \ assumes the form \ (\ref{1.2}),
\ with a suitably normalized \ Frobenius solution \ $\phi(\cdot,\lambda).$

We now give the formulas for all cases of Frobenius solutions which can occur
at $x=0$ under the assumptions (\ref{2.1})-(\ref{2.2}); these are the solutions
which were utilized in \cite{MN,FL}. The indicial equation for the R.S.P.
$x=0$ for the Sturm-Liouville equation(\ref{1.1}) with potential (\ref{2.1}),

\begin{equation}
-y^{\prime\prime}(x)+\left(  \frac{q_{0}}{x^{2}}+\frac{q_{1}}{x^{{}}}%
+\sum_{n=0}^{\infty}\,q_{n+2}x^{n}\text{ }\right)  \text{\ }y(x)=\lambda
y(x),\text{ \ \ \ \ \ }x\in(0,\infty).\label{2.4}%
\end{equation}
is%
\begin{equation}
r^{2}-r-q_{0}=r^{2}-r-\left(  \nu^{2}-\frac{1}{4}\right)  =\left(  r-\left(
\frac{1}{2}+\nu\right)  \right)  \cdot\left(  r-\left(  \frac{1}{2}%
-\nu\right)  \right)  \label{2.5}%
\end{equation}
where we have set
\[
q_{0}=\nu^{2}-\frac{1}{4},\nu\geq0
\]
for convenience. This gives rise to the following cases of Frobenius solutions:

\textbf{Case I:} $-\frac{1}{4}<q_{0}<\infty$, $\ \ q_{0}=\nu^{2}-\frac{1}%
{4}\neq\frac{M^{2}-1}{4},\ M=1,2,\dots$.(This is Case I\thinspace A in
\cite{MN}). \ In this case,
\begin{align}
y_{1}(x,\lambda)  & =x^{\frac{1}{2}+\nu}\left(  1+\sum_{n=1}^{\infty}%
\,a_{n}(\lambda)x^{n}\right)  ,\label{2.6}\\
y_{2}(x,\lambda)  & =x^{\frac{1}{2}-\nu}\left(  1+\sum_{n=1}^{\infty}%
\,b_{n}(\lambda)x^{n}\right)\label{2.7}  ,
\end{align}
where $a_{n}(\lambda),b_{n}(\lambda)$ are polynomials in $\lambda$ of degree
$\left[  \frac{n}{2}\right]  $, and
\begin{equation}
W_{x}\left(  y_{1}(\cdot;\lambda),y_{2}(\cdot;\lambda)\right)  =-2\nu
.\label{2.8}%
\end{equation}

\textbf{Case II A:} $q_{0}=\frac{M^{2}-1}{4},\ M$ odd: $M=2\ell+1,\,\ell
=0,1,\dots,$that is, $\ q_{0}=\ell(\ell+1).$

In \cite{MN} this is Case IC for M odd and it includes Case II (for
$\ell=0).$  In this case,%
\begin{align}
y_{1}(x,\lambda)  & =x^{\ell+1}\left(  1+\sum_{n=1}^{\infty
}\,a_{n}(\lambda)x^{n}\right)  ,\label{2.9}\\
y_{2}(x,\lambda)  & =K_{\ell}(\lambda)\,y_{1}(x;\lambda)\ln x+x^{-\ell}\left(
1+\sum_{n=1}^{\infty}\,d_{n}(\lambda)x^{n}\right)\label{2.10}  ,
\end{align}
where $a_{n}(\lambda),\,d_{n}(\lambda)$ are polynomials in $\lambda$ of degree
$\left[  \frac{n}{2}\right]  $, $K_{\ell}(\lambda)$ is a polynomial of degree
$\ell$, and
\begin{equation}
W_{x}\left(  y_{1}(\cdot,\lambda),y_{2}(\cdot,\lambda)\right)  =-(2\ell
+1).\label{2.11}%
\end{equation}

\textbf{Case II B:} $q_{0}=\frac{M^{2}-1}{4},\ M$ even: $M=2N,\,N=0,1,\dots
$, that is, $\ q_{0}=N^{2}-\frac{1}{4}.$\\[2mm]In \cite{MN} this is Case IC for
M even and it includes Case IB  (for $N=0$). \ In this case,%
\begin{align}
y_{1}(x,\lambda) &  =x^{\frac{1}{2}+N}\left(  1+\sum_{n=1}^{\infty}%
\,a_{n}(\lambda)x^{n}\right)  ,\label{2.12}\\
y_{2}(x,\lambda) &  =y_{1}(x,\lambda)\ln x+\sum_{n=1}^{\infty}\,d_{n}%
(\lambda)x^{\frac{1}{2}+n},\text{ \ if \ }N=0 \nonumber\\
y_{2}(x,\lambda) &  =K_{N}(\lambda)\,y_{1}(x;\lambda)\,\ln x\ +\ x^{\frac
{1}{2}-N}\left(  1+\sum_{n=1}^{\infty}\,d_{n}(\lambda)x^{n}\right)  ,\text{
\ if \ }N\geq1\label{2.13}%
\end{align}
%\begin{align}
%y_{1}(x,\lambda)  & =x^{\frac{1}{2}+N}\left(  1+\sum_{n=1}%
%^{\infty}\,a_{n}(\lambda)x^{n}\right)  ,\label{2.12}\\
%y_{2}(x,\lambda)  & =K_{N}(\lambda)\,y_{1}(x;\lambda)\,\ln x\ +\ x^{\frac
%{1}{2}-N}\left(  1+\sum_{n=1}^{\infty}\,d_{n}(\lambda)x^{n}\right)\label{2.13}  ,
%\end{align}
where $a_{n}(\lambda),\,d_{n}(\lambda)$ are polynomials in $\lambda$ of degree
$\left[  \frac{n}{2}\right]  $, $K_{N}(\lambda)$ is a polynomial of degree
$N,$ and
\begin{equation}
W_{x}\left(  y_{1}(\cdot,\lambda),y_{2}(\cdot,\lambda)\right)  =\left\{  \!\!%
\begin{array}
[c]{cll}%
-2N\text{ \ \ \ } & if & N\geq1,\\[1mm]%
1\text{ \ } & if & N=0.
\end{array}
\right.  \label{2.14}%
\end{equation}

In each of the above cases the first Frobenius solution $y_{1}(\cdot,\lambda)$
is the principal solution at $x=0$ for all $\lambda\in(-\infty,\infty)$. The
Frobenius solutions as normalized above satisfy the following properties (see
\cite[Theorem 2.1]{MN}):

\begin{itemize}
\item[{\bf (i)}] $y_{1}(x,\cdot),\,y_{2}(x,\cdot)$ and their derivatives are entire
functions for each $x\in(0,\infty)$ and satisfy for all $\lambda\in\Complex,\,x\in(0,\infty)$ the relations
\[
y_{i}\left(  x,\overline{\lambda}\right)  =\overline{y_{i}(x,\lambda)},\text{
\ \ \ }y_{i}^{\prime}\left(  x,\overline{\lambda}\right)  =\overline
{y_{i}^{\prime}(x,\lambda)},\text{ \ \ \ }i=1,2.
\]

\item[{\bf (ii)}] $y_{1}(\cdot,\lambda)\in L_{2}(0,x_{0})$ for $0<x_{0}<\infty$ and
for all $\lambda\in \Complex$.\vspace*{1mm}

\item[{\bf (iii)}] $W_{x}\left(  y_{1}(\cdot,\lambda),y_{2}(\cdot,\lambda)\right)
=C\neq0$ where $C\in R$, independent of $\lambda$.\vspace*{1mm}
\end{itemize}
\noindent It follows that a fundamental system of solutions \ near \ $x=0$ \ which is
entire in \ $\lambda$ \ and satisfies property \ {\bf (i)}, together with the
normalization

\begin{equation}
W_{x}\left(  \phi(\cdot,\lambda),\theta(\cdot,\lambda)\right)  =1\text{
\ }for\,\ all\text{ \ \ }\lambda\in \Complex,\label{2.15}%
\end{equation}
can be selected by taking%
\begin{equation}
\phi(x,\lambda):=y_{1}(x,\lambda),\text{ \ \ \ }\theta(x,\lambda
):=y_{2}(x,\lambda)/C,\label{2.16}%
\end{equation}
where \ C \ is the real constant in \ (\ref{2.8}), (\ref{2.11}), \ or (\ref{2.14}).
In the \textbf{\ LP} \ cases at $x=0$ \ only \ $\phi(\cdot,\lambda)$
\ satsifies \ the property \ {\bf (ii)} of square integrability \ near \ $x=0$, so
in all\textbf{ LP} cases it is the first Frobenius solution \ which is used
to write the eigenfunction expansion in the form
(\ref{1.2}).

The \textbf{LC}  cases at $x=0$ \ are \ Case I \ with \ $q_0\in(-\frac{1}%
{4},0)\cup(0,\frac{3}{4}),$ Case IIA with \ $q_0=0 (\ell = 0)$, \ and \ Case IIB with
\ $q_0=\frac{1}{4} (N=0)$.  In this paper we limit our consideration of {\bf LC} boundary 
conditions at $x=0$ to the Friedrichs {\bf LC} boundary condition (see (\ref{5.2}) below); in
this case it is the first Frobenius solution $\phi$ (the principal solution) which is selected
and used in the eigenfunction expansion (\ref{1.2}). 

% In these cases \ we can \ associate \ with equation
%\ref{1.1} the \ \textbf{LC} \ boundary condition 
%\begin{equation}
%W_{0}\left(  y,\theta(x,\lambda_{0})\right)  cos(\alpha)-W_{0}\left(
%y,\phi(x,\lambda_{0})\right)  sin(\alpha)=0,\text{ \ \ }\alpha\in\lbrack
%0,\pi).\label{2.17}%
%\end{equation}
%In these cases the solution \ $\phi=\phi_{\alpha}(\cdot,\lambda)$ \ to be used
%in the eigenfunction expansion \ \ref{1.2} \ for the \ problem \ \ref{1.1}%
%\ref{2.17} \ on \ $(0,\infty)$ \ can be defined by the following
% "\textit{LC initial conditions" \ }$\ $(see \cite{FULT}) \ %
%
%\begin{equation}
%\lim_{x\rightarrow0}\left(
%\begin{array}
%[c]{l}%
%W_{x}(\phi_{\alpha}(\cdot,\lambda),\phi(\cdot,0))\\
%W_{x}(\phi_{\alpha}(\cdot,\lambda),\theta(\cdot,0))
%\end{array}
%\right)  =\left(
%\begin{array}
%[c]{l}%
%\cos\alpha\\
%\sin\alpha
%\end{array}
%\right) \label{2.18}
%\end{equation}X
%so as to generate a solution satisfying \ property  (i) \ and \ the \textbf{LC}
%\ boundary condition \ref{2.17} \ for all $\lambda\in C.$ \ This yields,
%making use of the wronskian limit relations given in \cite[Theorem 2.1]{MN},%
%
%\begin{equation}
%\phi_{\alpha}(x,\lambda)=-\cos\alpha\text{ }\theta(x,\lambda)+\sin\alpha
%\phi(x,\lambda)\label{2.19},
%\end{equation}
%which is the solution to be used in the eigenfunction expansion (\ref{1.2}) for the
%problem (\ref{1.1}),(\ref{2.17}).
%where \{$\phi(x,\lambda),\theta(x,\lambda)$\} are defined \ (in the\textbf{LC}
%cases) \ as in \ref{2.16}. \ The Friedrichs LC boundary condition at \ x=0
%\ corresponds to \ $\alpha=\frac{\pi}{2}$ \ in \ref{2.17} and \ref{2.19}.

\section{Preliminaries}
\setcounter{equation}{0} In this section we collect together some useful
results which relate solutions of the Sturm-Liouville equation (\ref{1.1}) to
solutions of the companion first order system (\ref{1.6}). The proofs of
these results (though sometimes tedious) require only straightforward
algebraic manipulation making use of these two equations, and no special
assumptions on the potential $q(x)$; so we omit the proofs.

\vskip6pt \noindent\textbf{1.} If $y$ is any solution of the SL-equation, then
$(\,(y^{\prime})^{2},-2yy,y^{2})^{T}$ is a solution of the first order system
(\ref{1.6}).

\vskip6pt \noindent\textbf{2}. If we let a fundamental system of the
Sturm-Liouville equation (\ref{2.4}) be defined by the initial conditions at any $x_{0} > 0$
\begin{equation}
\left[
\begin{array}
[c]{ll}%
u(x_{0},\lambda) & v(x_{0},\lambda)\\
u^{\prime}(x_{0},\lambda) & v^{\prime}(x_{0},\lambda)
\end{array}
\right]  =\left[
\begin{array}
[c]{cc}%
1 & 0\\
0 & 1
\end{array}
\right]  ,\label{3.1}%
\end{equation}
then a corresponding fundamental system of solutions of equation (\ref{1.6})
is
\begin{equation}
U=\left[  U_{1},U_{2},U_{3}\right]  =\left[
\begin{array}
[c]{ccc}%
(u^{\prime})^{2} & u^{\prime}v^{\prime} & (v^{\prime})^{2}\\
-2uu^{\prime} & -[u^{\prime}v+uv^{\prime}] & -2vv^{\prime}\\
u^{2} & uv & v^{2}%
\end{array}
\right] .\label{3.2}%
\end{equation}

\vskip6pt \noindent\textbf{3. }$\ \ $If \  \{$\phi(x,\lambda),\theta
(x,\lambda)$\} are the Frobenius solutions defined in \ Section 2 (in all the
cases) \ and \ normalized \ by \ (\ref{2.16}) \ so as to ensure that
\ $W_{x}\left(  \phi(\cdot,\lambda),\theta(\cdot,\lambda)\right)  =1$, then a
corresponding fundamental system of solutions of equation (\ref{1.6}) \ is%

\begin{equation}
U=\left[  U_{1},U_{2},U_{3}\right]  =\left[
\begin{array}
[c]{ccc}%
(\theta^{\prime})^{2} & \theta^{\prime}\phi^{\prime} & (\phi^{\prime})^{2}\\
-2\theta\theta^{\prime} & -[\theta^{\prime}\phi+\theta\phi^{\prime}] &
-2\phi\phi^{\prime}\\
\theta^{2} & \theta\phi & \phi^{2}%
\end{array}
\right]  \label{3.3}%
\end{equation}

\vskip6pt \noindent\textbf{4. \ } An indefinite inner product on the solution
space \ of equation \ (\ref{1.6}) \ may be defined by
\begin{equation}
\langle U_{1},U_{2}\rangle:=2(P_{1}R_{2}+P_{2}R_{1})-Q_{1}Q_{2}%
=const,\mbox{ independent of x }\in\lbrack0,\infty)\label{3.4}%
\end{equation}
where $U_{k}=(P_{k},Q_{k},R_{k}),k=1,2$.

\vskip6pt \noindent\textbf{5.} \ For any solution $U$ = (P,Q,R)$^{T}$ of
equation (\ref{1.6}),%
\begin{equation}
\frac{d}{dx}\langle U,U\rangle=\frac{d}{dx}[4PR-Q^{2}]=0,\nonumber
\end{equation}
i.e.
\begin{equation}
4PR-Q^{2}=const,\mbox{ independent of x }\in\lbrack0,\infty
)\label{3.5}%
\end{equation}

\vskip6pt \noindent\textbf{6.} \ If $U_{1}$ and $U_{2}$ are any two solutions
of equation (\ref{1.6}) represented in the form,
\begin{equation}
U_{j}=\left[
\begin{array}
[c]{c}%
P_{j}\\
Q_{j}\\
R_{j}%
\end{array}
\right]  =a_{j}\left[
\begin{array}
[c]{c}%
(\theta^{\prime})^{2}\\
-2\theta\theta^{\prime}\\
\theta^{2}%
\end{array}
\right]  +b_{j}\left[
\begin{array}
[c]{c}%
\theta^{\prime}\phi^{\prime}\\
-[\theta\phi^{\prime}+\theta^{\prime}\phi]\\
\theta\phi
\end{array}
\right]  +c_{j}\left[
\begin{array}
[c]{c}%
(\phi^{\prime})^{2}\\
-2\phi\phi^{\prime}\\
\phi^{2}%
\end{array}
\right]  \label{3.6}
\end{equation}
in terms of the Frobenius solutions \{${\phi(\cdot,\lambda),\theta
(\cdot,\lambda)}$\} of section 2, then we have%

\begin{equation}
\langle U_{1},U_{2}\rangle:=2(P_{1}R_{2}+P_{2}R_{1})-Q_{1}Q_{2} = 2(a_{1}%
c_{2}+c_{1}a_{2})-b_{1}b_{2}\label{3.7}%
\end{equation}

In particular,
\begin{equation}
\langle U_{1},U_{1}\rangle:= 4P_{1}R_{1}-Q_{1}^{2}= 4a_{1}c_{1}-b_{1}%
^{2}.\label{3.8}%
\end{equation}

\vskip6pt \noindent\textbf{7.} Similarly, if these same solutions, $U_{1}$ and
$U_{2}$, of equation (\ref{1.6}) are represented in the form,
\begin{equation}
U_{j}=\left[
\begin{array}
[c]{c}%
\tilde{P}_{j}\\
\tilde{Q}_{j}\\
\tilde{R}_{j}%
\end{array}
\right]  =\tilde{a_j}\left[
\begin{array}
[c]{c}%
(u^{\prime})^{2}\\
-2%
%TCIMACRO{\unit{u}}%
%BeginExpansion
\operatorname{u}%
%EndExpansion
u^{\prime}\\
u^{2}%
\end{array}
\right]  +\tilde{b_j}\left[
\begin{array}
[c]{c}%
u^{\prime}v^{\prime}\\
-[uv^{\prime}+u^{\prime}v]\\
uv
\end{array}
\right]  +\tilde{c_j}\left[
\begin{array}
[c]{c}%
(v^{\prime})^{2}\\
-2vv^{\prime}\\
v^{2}%
\end{array}
\right]  \label{3.9}%
\end{equation}
in terms of the solutions defined in (\ref{3.1})-(\ref{3.2}) we have
\begin{equation}
\langle U_{1},U_{2}\rangle:=2(\tilde{P}_{1}\tilde{R}_{2}+\tilde{P}_{2}\tilde{R}_{1})-\tilde{Q}_{1}\tilde{Q}_{2}=2(\tilde{a}_{1}%
\tilde{c}_{2}+\tilde{c}_{1}\tilde{a}_{2})-\tilde{b}_{1}\tilde{b}_{2}.\label{3.10}%
\end{equation}
In particular,
\begin{equation}
\langle U_{1},U_{1}\rangle:= 4\tilde{P}_{1}\tilde{R}_{1}-\tilde{Q}_{1}^{2}= 4\tilde{a}_{1}\tilde{c}_{1}-\tilde{b}_{1}^{2}.\label{3.11}%
\end{equation}
It follows from (\ref{3.7}), (\ref{3.10})  and  (\ref{3.8}), (\ref{3.11}) that we must also have %
\begin{align}
2(a_{1}c_{2}+c_{1}a_{2})-b_{1}b_{2}  & =2(\tilde{a}_{1}\tilde{c}_{2}+\tilde{c}_{1}\tilde{a}_{2})-\tilde{b}_{1}\tilde{b}_{2}, \text{ \ \ and}
4a_{1}c_{1}-b_{1}^{2}  & =4\tilde{a}_{1}\tilde{c}_{1}-\tilde{b}_{1}^{2}.
\end{align}
\vskip6pt \noindent\textbf{8. \ \ }If $y$ is any \ solution \ of \ the \ SL
\ equation (\ref{1.1}) and $U$=(P,Q,R)$^{T}$ is any solution of the companion
system (\ref{1.6}) \ then
\begin{equation}
\frac{d}{dx}[Py^{2}+Qyy^{\prime}+R(y^{\prime})^{2}]=0,\label{3.14}%
\end{equation}
i.e.,%
\begin{equation}
P(x,\lambda)y^{2}(x,\lambda)+Q(x,\lambda)y(x,\lambda)y^{\prime}(x,\lambda
)+R(x,\lambda)(y^{\prime}(x,\lambda)^{2}=\mbox{ constant, independent of \
x.}\label{3.15}%
\end{equation}

% This is section4
\section{A Spectral Density  Function Characterization of Al-Naggar and Pearson}
\setcounter{equation}{0}
In this section we consider the Sturm-Liouville problem
\begin{equation}
-y^{\prime\prime}(x)+\left(  \frac{q_{0}}{x^{2}}+\frac{q_{1}}{x^{{}}}%
+\sum_{n=0}^{\infty}\,q_{n+2}x^{n}\text{ }\right)  \text{\ }y(x)=\lambda
y(x),\text{ \ \ \ \ \ }x\in[A,\infty), \text{ \ \ }A > 0,\label{4.1}%
\end{equation}
\begin{equation}
y(A) = 0.\label{4.2}%
\end{equation}
We shall assume that {\bf Assumption 1, Case I,} holds, so that the above potential $q$ is continuous in $(0,\infty)$ and also has a continuous derivative.  In addition we make the following assumption:
\begin{center}
{\bf Assumption 3:} $\mathbf{Near\ x=\infty:}$
\end{center}
For $x_0 > 0$ we have either
\begin{equation}
q \in L_1(x_0,\infty) \label{4.3}
\end{equation}
or
\begin{equation}
q' \in L_1(x_0,\infty), \qquad q \in AC_{loc}[x_0,\infty),  \qquad  {\text and}\lim_{x\rightarrow\infty}q(x)=0. \label{4.4}
\end{equation}

Under the assumption (\ref{4.3}) or the assumption (\ref{4.4}) it was established in Fulton, Pearson and Pruess \cite[Thm1 and Cor 2]{FPP2} that the
initial value problem (\ref{1.6})-(\ref{1.7}) has a unique solution for all $\lambda \in (0,\infty)$.  Henceforth we denote
this unique solution by

\begin{equation}
U_1(x,\lambda)=\left[
\begin{array}
[c]{c}%
P_1(x,\lambda)\\
Q_1(x,\lambda)\\
R_1(x,\lambda)
\end{array}
\right], \label{4.5}
\end{equation}
that is, $U_1$ is the unique solution (under {\bf Assumption 3}) for which
\[
\lim_{x\rightarrow\infty} U_1(x,\lambda)
=  \left(
\begin{array}
[c]{l}
P_1(\infty,\lambda)\\
Q_1(\infty,\lambda)\\
R_1(\infty,\lambda)
\end{array}
\right)
=\left(
\begin{array}
[c]{l}
\sqrt{\lambda}\\
0\\
\frac{1}{\sqrt{\lambda}}%
\end{array}
\right)  ,\text{ \ \ \ for \ }\lambda\in(0,\infty).
\]

The assumptions (\ref{4.3}) and (\ref{4.4}) were used in \cite{FPP2,FPP4} to obtain spectral density function characterizations of the type (\ref{1.5}) when the left endpoint is regular. Al-Naggar and Pearson \cite{AlN1,AlN2} also obtained a spectral density function characterization of this type when the left endpoint is regular, using a different approach.  
Their approach was based on the determination of intervals of a.c. spectrum by locating those intervals of the real $\lambda$-axis where subordinate solutions do not exist. The purpose of this section is to present their analysis as it applies to the problem (\ref{4.1})-(\ref{4.2}) and under the additional {\bf Assumption 3} on the half line $[A,\infty)$, $A > 0$, and show that it also guarantees a.c. spectrum for $\lambda \in (0,\infty)$ and yields a spectral density function formula of the type (\ref{1.5}) (see (\ref{4.30}) below and \cite[Cor 4]{FPP2}); this route to the spectral density function characterization represents an alternative to the analysis in \cite{FPP2,FPP4}.  As we shall see, the approach of Al-Naggar and Pearson has a major advantage in that it  extends nicely to obtain a corresponding spectral density function characterization of the type (\ref{1.5}) for the doubly singular equation (\ref{4.1}) on $(0,\infty)$.  This will be done in the next section.

In \cite{AlN1,AlN2} the analysis  is focused on the third order ordinary differential equation (see (\ref{4.13}) below) satisfied by the third component of a solution of Appell's first order system; here, we modify the approach slightly so as to focus attention on the system (\ref{1.6}), so that we can properly exploit the results from \cite{FPP2} on uniqueness of the above solution $U_1$ satisfying the
initial condition (\ref{1.7}).

Letting $\{u(\cdot,\lambda), v(\cdot,\lambda)\}$ be the fundamental system of solutions of (\ref{4.1}) defined by the initial conditions (\ref{3.1}) at $x_0=A$, the Titchmarsh-Weyl m-function associated with the problem (\ref{4.1})-(\ref{4.2}) is defined by
\begin{equation}
 \Psi_{A}(\cdot,\lambda) := u(\cdot,\lambda) + m_A(\lambda)v(\cdot,\lambda) \in L_2(A,\infty),  \text { \ \ } \text{for all $Im(\lambda) \ne 0$}.\label{4.6}
\end{equation}
Then, as is well known, this m-function is a Nevanlinna function 
%(also {\it Pick-Nevanlinna} or {\it Herglotz})
and therefore admits the representation
\begin{equation}
m_A(z)=\alpha+\beta z+\int_{-\infty}^{\infty} \left(  \sfrac{1}{t-z}-\sfrac{t}{1+t^{2}}\right)   d\rho_A(t),\text{ \ \ } \alpha \in (-\infty,\infty), \text{ \ \ } \beta \ge 0,\label{4.7}%
\end{equation}
where  the inversion integral for $\rho_A$ in terms of $m$ is the Titchmarsh-Kodaira formula
\begin{equation}\label{4.8}
\rho_A(\lambda)=\lim_{\varepsilon\downarrow 0} \sfrac {1}{\pi} \int_0^{\lambda} Im[m_A(t+{\rm i}\varepsilon)] dt.
\end{equation}
Since the {\bf Assumption 3} ensures a.c. spectrum for $\lambda \in (0,\infty)$ we may differentiate (\ref{4.8}) to obtain the spectral density function as
\begin{equation}
f_{A}(\lambda) := \rho^{\prime}_{A}(\lambda) = \lim_{\epsilon\downarrow0}\frac{1}{\pi} Im[m_{A}(\lambda+i\epsilon)] \label{4.9}
\end{equation}  

We now proceed to transform (\ref{4.9}) into the form (\ref{1.5}); after some analysis this yields (\ref{4.30}) in Theorem 1 below. In the sequel it will be helpful to make use of the fundamental system of solutions of Appell's system
(\ref{1.6}) given in (\ref{3.2}) in terms of the solutions $\{u(\cdot,\lambda), v(\cdot,\lambda)\}$ of equation (\ref{4.1}) fixed by the initial conditions (\ref{3.1}) with $x_0=A$

\vskip 6pt \noindent\textsc{Lemma 1.}  Assume the potential $q$ satisfies {\bf Assumption 3}.\newline
{\bf (i)} For all $\lambda \in (0,\infty)$ let $U_1$ be
be the unique solution defined at $x=\infty$ in (\ref{4.5}). Then using the indefinite inner product on the solution space of (\ref{1.6}) defined in (\ref{3.7}) we have
\begin{equation}\label{4.10}
\langle U_1, U_1 \rangle= 4 P_1(x, \lambda) R_1(x, \lambda) - (Q_1(x,\lambda))^{2} = 4.
\end{equation}
{\bf (ii)} If the solution $U_1$ is represented in the form (\ref{3.9}), say
\begin{equation}
U_1(x,\lambda)=\left[
\begin{array}
[c]{c}%
P_1(x,\lambda)\\
Q_1(x,\lambda)\\
R_1(x,\lambda)
\end{array}
\right]  = \tilde a\left[
\begin{array}
[c]{c}%
(u^{\prime})^{2}\\
-2uu^{\prime}\\
u^{2}%
\end{array}
\right]  + \tilde b \left[
\begin{array}
[c]{c}%
u^{\prime}v^{\prime}\\
-[uv^{\prime}+ u^{\prime}v]\\
uv
\end{array}
\right]  + \tilde c\left[
\begin{array}
[c]{c}%
(v^{\prime})^{2}\\
-2vv^{\prime}\\
v^{2}%
\end{array}
\right]  \label{4.11}
\end{equation}
then 
\begin{equation}
 4 \tilde a \tilde b - (\tilde c)^2 = 4. \label{4.12}
\end{equation}
%{\bf (iii)} If the solution $U_1$ is represented in the form (\ref{3.6}), say
%\begin{equation}
%U_1(x,\lambda)=\left[
%\begin{array}
%[c]{c}%
%P_1(x,\lambda)\\
%Q_1(x,\lambda)\\
%R_1(x,\lambda)
%\end{array}
%\right]  =a^{*}\left[
%\begin{array}
%[c]{c}%
%(\theta^{\prime})^{2}\\
%-2\theta\theta^{\prime}\\
%\theta^{2}%
%\end{array}
%\right]  +b^{*}\left[
%\begin{array}
%[c]{c}%
%\theta^{\prime}\phi^{\prime}\\
%-[\theta\phi^{\prime}+\theta^{\prime}\phi]\\
%\theta\phi
%\end{array}
%\right]  +c^{*}\left[
%\begin{array}
%[c]{c}%
%(\phi^{\prime})^{2}\\
%-2\phi\phi^{\prime}\\
%\phi^{2}%
%\end{array}
%\right]  \label{4.13}
%\end{equation}
%then
%\begin{equation}
%4 a^{*} c^{*} - (b^{*})^{2} = 4. \label{4.14}
%\end{equation}
\noindent\textsc{Proof.}{\bf (i)} To get the constant 4 observe from (\ref{3.5}) and
the initial condition (\ref{1.7}) that
\[ 4 P_1(x, \lambda) R_1(x, \lambda) - (Q_1(x,\lambda))^{2} = \lim_{x\to\infty} 4 P_1(x, \lambda) R_1(x, \lambda) - (Q_1(x,\lambda))^{2} = 4\sqrt{\lambda}\cdot \big(\sfrac{1}{\sqrt{\lambda}}\big)=4.\]
{\bf (ii)} The conversion of the indefinite inner product in terms of the 
coefficients $\{ \tilde a, \tilde b, \tilde c\}$ is the property (\ref{3.7}) and (\ref{3.8}); this is readily proved
by substitution of the components of $U_1$ into the inner product formula and use of the wronskian
relation  $W_x(u(\cdot,\lambda),v(\cdot,\lambda))=1$.  \quad \rule{2.2mm}{3.2mm}

It can be shown that the third component, $R(x,\lambda)$, of a solution of Appell's system (\ref{1.6}) satisfies the third order
equation  of \cite[p. 6584, Equa. 5]{AlN2},
\begin{equation}
\frac{d^{3}R}{dx^{3}} + 4(\lambda-
q(x)) \frac{dR}{dx} - 2 q^{\prime}(x) R = 0. \label{4.13}
\end{equation}

{\bf Remark} In Appell's paper this third order equation is \cite[p. 213, Equa (5)]{APPELL}.  
\vskip 6pt
\noindent  
Also, a general solution of (\ref{4.13}) can be obtained from a suitably normalized solution $R_{1}(x, \lambda)$ as in \cite[p. 6587, Lemma~2]{AlN2}. Here we generalize this technique to obtain from the given solution $U_{1}$, two other linearly independent solutions. This is the content of the following lemma.

\vskip 6pt \noindent\textsc{Lemma 2.} Let $U_{1} = (P_{1}, Q_{1}, R_{1})^T$ be the unique solution of (\ref{1.6}) defined at $x=\infty$ by the initial condition (\ref{1.7}).
% Then the identity (\ref{4.10}) holds for all $x \in (0,\infty)$; moreover, 
Then we can write the  
general solution of (\ref{1.6}) in the form
\begin{equation}
\label{4.14}U = \left[
\begin{array}
[c]{c}%
P\\
Q\\
R
\end{array}
\right]  = \beta_1 U_{1} + \beta_2 U_{2} +\beta_3 U_{3},
\end{equation}
where
\begin{equation}
\label{4.15}U_{2}(x, \lambda) = \left[
\begin{array}
[c]{c}%
P_{1} \cos2 \gamma+ (Q_{1} / R_{1}) \sin2 \gamma- (2 / R_{1}) \cos2 \gamma\\
Q_{1} \cos2 \gamma+ 2 \sin2 \gamma\\
R_{1} \cos2 \gamma
\end{array}
\right] ,
\end{equation}
\begin{equation}
\label{4.16}U_{3}(x, \lambda) = \left[
\begin{array}
[c]{c}%
P_{1} \sin2 \gamma- (Q_{1} / R_{1}) \cos2 \gamma- (2 / R_{1}) \sin2 \gamma\\
Q_{1} \sin2 \gamma- 2 \cos2 \gamma\\
R_{1} \sin2 \gamma
\end{array}
\right] ,
\end{equation}
and
\[
\gamma(x) = \int_{x_{0}}^{x} 1 / R_{1}(t, \lambda) \, dt
\]
for some $x_{0} > 0$.

\bigskip

\noindent\textsc{Proof.} The fact that
\[
R(x, \lambda) := R_{1}(x, \lambda) [\beta_1 + \beta_2 \cos2 \gamma+ \beta_3 \sin2 \gamma]
\]
is a general solution of (\ref{4.13}) follows as in \cite[Lemma 2, p.
6587]{AlN2}. Since (\ref{4.13}) is satisfied by the third component of any
solution of (\ref{1.6}), we can generate a general solution for (\ref{1.6})
by computing $Q(x, \lambda) = -dR/dx$ and then $P(x, \lambda) = [-dQ/dx + 2
(\lambda- q) R] / 2$ and expressing these solutions of (\ref{1.6}) in terms
of $P_{1}$, $Q_{1}$, $R_{1}$. This gives the result (\ref{4.14}).
Alternatively, a direct substitution of $U_{2}$ and $U_{3}$ into (\ref{1.6})
and use of the formulas for $P_{1}^{\prime}$, $Q_{1}^{\prime}$ and
$R_{1}^{\prime}$ will verify that $U_{2}$ and $U_{3}$ are solutions of
(\ref{1.6}). \quad\rule{2.2mm}{3.2mm}

\vskip 6pt \noindent\textsc{Corollary.} {\bf (i)} The general solution (\ref{4.14}) in
Lemma~2 satisfies
\begin{equation}
\label{4.17}\langle U, U \rangle= 4 PR - Q^{2} = 4 ((\beta_1)^{2} - (\beta_2)^{2} - (\beta_3)^{2}).
\end{equation}
Similarly, if
\[
\tilde U = (\tilde P, \tilde Q, \tilde R) = \tilde \beta_1  U_{1} + \tilde \beta_2 U_{2} +
\tilde \beta_3 U_{3},
\]
we have for the inner product that
\begin{align}
\langle U, \tilde U \rangle & = 2 (P \tilde R + R \tilde P) - Q \tilde
Q\label{4.18}\\
& = 4 (\beta_1 \tilde \beta_1 - \beta_2 \tilde \beta_2 - \beta_3 \tilde \beta_3).\nonumber
\end{align}
{\bf (ii)} The solutions $\{U_1,U_2,U_3\}$ in Lemma 1 are mutually orthogonal with respect to the
indefinite inner product defined in (\ref{3.4}).

\noindent\textsc{Proof.} {\bf (i)} For (\ref{4.17}) calculate $4PR - Q^{2}$ using
(\ref{4.14}) and the normalization (\ref{4.10}). This simplifies to
\[
\begin{array}
[c]{rcl}
4(\beta_1 + \beta_2 \cos2 \gamma+ \beta_3 \sin2 \gamma)^{2} & - & 8 (\beta_2 \cos2 \gamma+ \beta_3 \sin2
\gamma) (\beta_1 + \beta_2 \cos2 \gamma+ \beta_3 \sin2 \gamma)\\
& - & 4 ((\beta_2)^{2} \sin^{2} 2 \gamma+ (\beta_3)^{2} \cos^{2} 2 \gamma- 2 \beta_2 \beta_3 \sin2
\gamma\cos2 \gamma)\\
& = & 4 ((\beta_1)^{2} - (\beta_2)^{2} - (\beta_3)^{2}).
\end{array}
\]
the proof of (\ref{4.18}) is similar. \newline
\noindent
{\bf (ii)} Taking $U = U_{1}$, which statisfies the  normalization 
(\ref{4.10}), and $\tilde U = U_{2}$ in (\ref{4.15}) we have from (\ref{4.18})
that $\langle U_{1}, U_{2} \rangle= 4(1 \cdot0 - 0 \cdot1 - 0 \cdot0) = 0$.
Similarly, $\langle U_{1}, U_{3} \rangle$ and $\langle U_{2}, U_{3} \rangle$
are zero.
% Hence, the linearly independent solutions defined in Lemma~1 are
%mutually orthogonal in the sense of the indefinite inner product defined in (\ref{3.4}).
\quad\rule{2.2mm}{3.2mm}
\vspace{0.1in}

% The next lemma shows that the unique solution of the 
%initial value problem (\ref{1.6})-(\ref{1.7}) at  $x=\infty$ satisfies the normalization condition (\ref{5.18}) needed
%in Lemma 1.

\vspace{0.1in}
\noindent
{\bf Remark.} The third component,  $R_1$, of (\ref{4.11}) satisfies the third order equation (\ref{4.13}) and since
$R_1(x,\lambda) \rightarrow 1/\sqrt{\lambda} > 0$, it is in fact the same quadratic form as {\it Y(x,$\lambda$)}  
which was employed  in \cite[Lemma3 and Thm2]{AlN2}.
\vspace{0.1in}

The following lemma gives a limit relation 
involving the solution, $U_1$, and the solutions from Lemma 2 which are orthogonal to it.

\vskip 6pt \noindent\textsc{Lemma 3.} We assume the potential $q$ satisfies {\bf Assumption 3}. \newline
{\bf (i)} For the solution
$U_1$ of the initial value problem (\ref{1.6})-(\ref{1.7}), let $U_2$ and $U_3$ be the linearly independent solutions
of Lemma 1 which are generated by using $(P_1,Q_1,R_1)^{T}$ in (\ref{4.15}) and (\ref{4.16}). Then, for any linear
combination
\[
  U = \left(
\begin{array}
[c]{c}%
P\\
Q\\
R
\end{array} 
\right)
=\beta_2U_2 + \beta_3U_3
\]
we have  for all $x_0 > 0$,
\begin{equation}\label{4.19}
\lim_{N \to\infty} \frac{\int_{x_{0}}^{N} R(x,
\lambda) \, dx}{\int_{x_{0}}^{N} R_1(x, \lambda) \, dx} = 0.
\end{equation}
{\bf (ii)} The representation (\ref{4.11}) for $U_1$ has $\tilde a(\lambda) > 0$ and $\tilde c(\lambda) > 0$ for all $\lambda \in (0,\infty)$.
\vspace{0.1in}

\noindent\textsc{Proof} {\bf (i).} The solutions $U_{2}$ and $U_{3}$ in (\ref{4.15}) and 
(\ref{4.16}) have $R_{2} = R_{1}(x, \lambda) \cos( 2 \gamma(x))$ and  \newline
$R_{3} = R_{1}(x, \lambda) \sin( 2 \gamma(x))$ , where
\[
\gamma(x) = \int_{x_{0}}^{x} [1 / R_{1}(t, \lambda)] \, dt,
\]
for all $x > 0$. The fact that $R_1(x) > 0$ for all $x > 0$ follows from
(\ref{4.20}) below. Consequently, it suffices to prove that (\ref{4.19}) holds
for these two choices of $R$. For $R_{2}$ we have
\begin{align*}
\int_{x_{0}}^{N} R_{1} \cos2 \gamma\, dx  & = \int_{x_{0}}^{N} 0.5 (R_{1})^{2}
\frac{d}{dx} \sin2 \gamma\, dx\\
& = 0.5 (R_{1})^{2} \sin2 \gamma\big|_{x_{0}}^{N} + \int_{x_{0}}^{N} R_{1}
Q_{1} \sin2 \gamma\, dx.
\end{align*}
Since $Q_1(x) \to 0$ as $x \to\infty$ we can prove that
\[
\lim_{N \to\infty} \frac{\int_{x_{0}}^{N} R_{1} Q_{1} \sin2
\gamma\, dx}{\int_{x_{0}}^{N} R_{1} \, dx} = 0.
\]
Given $\epsilon> 0$, pick $x_{\epsilon}$ sufficiently large that $|Q_{1}(x)| <
\epsilon$ for $x \ge x_{\epsilon}$. Then we have
\[
\left|  \int_{x_{0}}^{N} R_{1} Q_{1} \sin2 \gamma\, dx \right|  \le M
\int_{x_{0}}^{x_{\epsilon}} R_{1} \, dx + \epsilon\int_{x_{\epsilon}}^{N}
R_{1} \, dx,
\]
where $M$ is a bound on $Q_1$ in $[x_{0}, x_{\epsilon}]$. Now pick $N$
sufficiently large that
\[
\frac{\displaystyle M \int_{x_{0}}^{x_{\epsilon}} R_{1} \, dx}%
{\displaystyle \int_{x_{0}}^{N} R_{1} \, dx} < \epsilon,
\]
which is possible since $\int_{x_{0}}^{\infty} R_{1} \, dx = \infty$. Then the
above quotient is less than $2 \epsilon$. For the boundary term we have
\begin{align*}
[R_1(N, \lambda)]^{2}  & = [R_1(x_{0}, \lambda)]^{2} + \int_{x_{0}}^{N} 2
R_1 R_1^{\prime}\, dx\\
& = [R_1(x_{0}, \lambda)]^{2} - \int_{x_{0}}^{N} 2 R_1 Q_1 \, dx,
\end{align*}
so
\[
\lim_{N \to\infty} \frac{\displaystyle [R_1(N, \lambda)]^{2}}%
{\int_{x_{0}}^{N} R_1 \, dx} = 0,
\]
by employing the same argument. The proof for $R_{3}$ is similar. \newline
\vspace{0.1in}
{\bf (ii)} In the representation (\ref{4.11}) for $U_1$ we have
\[
R_1(x, \lambda) = \tilde a u(x, \lambda)^{2} + \tilde b u(x, \lambda) v(x,
\lambda) + \tilde c v(x, \lambda)^{2}.
\]
To see that $\tilde a>0$ and $ \tilde c>0$ for all $\lambda\in(0, \infty)$ we first observe
from (\ref{4.12}) that $4 \tilde a \tilde c - (\tilde b)^{2} = 4$ requires $\tilde a$ and $\tilde c$ to be of the
same sign; and this must hold for all $\lambda \in (0,\infty)$ since they are continuous and cannot pass through zero (which would violate (\ref{4.12})). Factoring out $\tilde a$ we have
\[
R_1(x, \lambda) = \tilde a [u(x, \lambda)^{2} + (\tilde b/\tilde a) u(x, \lambda)
v(x, \lambda) + (\tilde c/\tilde a) v(x, \lambda)^{2}]
\]
with the coefficients of $u^{2}$ and $v^{2}$ positive. Therefore
$R_1$ admits a factorization of the form
\begin{align}
R_1(x, \lambda)  & = \tilde a[(u(x, \lambda) + \alpha v(x, \lambda))
(u(x, \lambda) + \bar\alpha v(x, \lambda))]\nonumber\\
& = \tilde a |u(x, \lambda) + \alpha v(x, \lambda)|^{2},\label{4.20}
\end{align}
where $\alpha= \alpha_{1} + i \alpha_{2}$ must satisfy (by (\ref{4.12}))
\begin{equation}
 \alpha_{1} = -\tilde b / (2\tilde a) \mbox{ and } \alpha_{2}^{2} = (\tilde c / \tilde a) - [(\tilde b)^{2} / (4(\tilde a)^{2})] = 1 / (\tilde a)^{2}. \label{4.21}
\end{equation}
Since $\tilde a$, $\tilde b$, and $\tilde c$ are real, we must have either $\alpha_{2} = 1/\tilde a$ or
$\alpha_{2} = -1 / \tilde a$; but either way the factorization remains the same with
$\alpha$ and $\bar\alpha$ switched. Since $|u-\alpha v|^{2}>0$, it
follows from the above factorization that $\tilde a(\lambda)>0$ for all $\lambda
\in(0, \infty)$; otherwise, $\lim_{x \to\infty} R_1(x, \lambda) \le0$
contradicting the fact that $R_1(\infty, \lambda) = 1 / \sqrt{\lambda} > 0$.
Hence, also $\tilde c(\lambda) > 0$ for all $\lambda\in(0, \infty)$. \quad
\rule{2.2mm}{3.2mm}

\vskip 6pt Since the definition of $\alpha_{2}$ must be $\pm1 / \tilde a(\lambda)$
from (\ref{4.21}) and since this indeterminacy is actually immaterial, we
choose to take $\alpha_{2} = 1 / \tilde a(\lambda)$, so that
\begin{equation}
\alpha(\lambda) := -\frac{\tilde b}{ 2 \tilde a} + i \;
\frac{ 1}{\tilde a}. \label{4.22}
\end{equation}
Using the fundamental system $\{u(\cdot,\lambda),v(\cdot,\lambda)\}$ we now define the complex-valued solution of (\ref{4.1}) for real $\lambda
\in(0, \infty)$,  
\begin{equation}
\label{4.23}\psi_A(x, \lambda) := u(x, \lambda) + \alpha(\lambda) v(x,\lambda).
\end{equation}
The key idea of Al-Naggar and Pearson, which enables identification of  a.c. spectrum, is embodied in
the following requirement:

\vskip 6pt \noindent{\bf Definition.} The general Sturm-Liouville equation (\ref{1.1}) satisfies 
{\em Condition A}, for a given real value of $\lambda$, if and only if 
there exists a complex-valued solution $y(x, \lambda)$ of (\ref{1.1})
for which
\begin{equation}    \label{4.24}
    \lim_{N \to \infty} \sfrac{\int_{x_0}^N  y(x, \lambda)^2  \, dx}
                        {\int_{x_0}^N |y(x, \lambda)|^2 \, dx} = 0 {\text \ \ \ for \ \  x_0 > 0}.
\end{equation}
\vspace{0.1in}
We now prove that equation (\ref{4.1}) satisfies {\em Condition A} for all $\lambda \in (0,\infty)$.
\vskip 6pt \noindent\textsc{Lemma 4.} Assume $q(x)$ in (\ref{4.1}) satisfies {\bf Assumptions 3}. Then for $x_{0} > 0$ and all $\lambda\in(0, \infty)$
\begin{equation}
\label{4.25}\lim_{N \to\infty} \frac{\int_{x_{0}}^{N} \psi_A(x,
\lambda)^{2} \, dx}{\int_{x_{0}}^{N} |\psi_A(x, \lambda)|^{2} \,
dx } = 0.
\end{equation}
%This means that $\psi_A(x, \lambda)$ as defined in (\ref{4.23}) satisfies the
%\textit{Condition~A} in (\ref{4.24}) for all $\lambda\in(0,\infty)$.
% it follows from {\bf Theorem 0} that the
%self-adjoint operator associated with the Sturm-Liouville problem  (\ref{4.1})-(\ref{4.2}) has
%a.c. spectrum in $(0,\infty)$.
\vspace{0.1in}

\noindent\textsc{Proof:} Since equation (\ref{4.13}) is satisfied by all
linear combinations of $v(\cdot, \lambda)^{2}$, $u(\cdot, \lambda)
v(\cdot, \lambda)$ and $u(\cdot, \lambda)^{2}$ (see \cite[Lemma~1, p.
6584]{AlN2}), it follows that $\psi_A(\cdot, \lambda)^{2}$ is a solution of
(\ref{4.13}) and since it is complex-valued, also that Re [$\psi_A(\cdot,
\lambda)^{2}$] and Im [$\psi_A(\cdot, \lambda)^{2}$] satisfy (\ref{4.13}).
Accordingly, it follows that there exist solutions $U_{2} = (P_{2}, Q_{2},
\mbox{Re}(\psi_A^{2}))^T$ and $U_{3} = (P_{3}, Q_{3}, \mbox{Im} (\psi_A^{2}))^T$ of the
first order system (\ref{1.6}), since any real solution $R$ of (\ref{4.13})
can be used to generate a solution of (\ref{1.6}) having $R$ as its third
component; e.g., let $Q^{\prime}= -R^{\prime}$ and $P = - Q^{\prime}/ 2 -
(\lambda- q)R$. From (\ref{4.20}) and (\ref{4.23}) we readily deduce the
following representations of $R_1(\cdot, \lambda)$, $Re[ \psi_A(\cdot,
\lambda)]^{2}$, and $Im[ \psi_A(\cdot, \lambda)]^{2}$ of the form (\ref{4.11}) (or, the third component of (\ref{4.11})):
\begin{align*}
R_1(x, \lambda)  & = \tilde a u^{2} + 2 \tilde a \alpha_{1} uv+ \tilde a (\alpha
_{1}^{2} + \alpha_{2}^{2}) v^{2}\\
\mbox{Re }[ \psi_A(x, \lambda)]^{2}  & =  u^{2} + 2 \alpha_{1} uv+
(\alpha_{1}^{2} - \alpha_{2}^{2}) v^{2}\\
\mbox{Im} [\psi_A(x, \lambda)]^{2}  & =  2 \alpha_{2} uv+ 2 \alpha_{1}\alpha_{2} v^{2}.\\
\end{align*}
It now follows from these formulas that the above solutions $U_{2}$ and $U_{3}$
associated with $Re[ \psi_A(\cdot, \lambda)]^{2}$ and $Im[ \psi_A(\cdot, \lambda
))]^{2}$ are orthogonal to $U_1$ in the sense of the inner product defined in (\ref{3.7}), i.e., we have
\[
\langle U_1, U_{2} \rangle= 2 [\tilde a (\alpha_{1}^{2} - \alpha_{2}^{2}) + \tilde a
(\alpha_{1}^{2} + \alpha_{2}^{2})] - 4 \tilde a \alpha_{1}^{2} = 0,
\]
and
\[
\langle U_1, U_{3} \rangle= 2 [ \tilde a (2 \alpha_{1} \alpha_{2}) +0] - 4 \tilde a
\alpha_{1} \alpha_{2} = 0.
\]
Hence it follows from Lemma~3(i) that for all $x_{0} > 0$
\[
\lim_{N \to\infty} \frac{ \int_{x_{0}}^{N} \mbox{Re } [\psi_A(x,
\lambda)]^{2} \, dx}{\int_{x_{0}}^{N} R_1(x, \lambda) \, dx} =0
\]
and
\[
\lim_{N \to\infty} \frac{ \int_{x_{0}}^{N} \mbox{Im } \psi_A(x,
\lambda)]^{2} \, dx}{\int_{x_{0}}^{N} R_1(x, \lambda) \, dx} =
0,
\]
from which (\ref{4.25}) follows. \quad\rule{2.2mm}{3.2mm}
\vspace{0.1in}

Assumming the Sturm-Liouville equation (\ref{1.1}) has a potential $q$ which is {\bf LP} at $x=\infty$ and regular at a finite left endpoint, Al-Naggar and Pearson obtain the following results in \cite[Lemma1 and Thm2]{AlN1} (where the fundamental system $\{u,v\}$ is defined by initial conditions at the left endpoint so that $v$ satisfies a general regular boundary condition and $ W_{x}(u(\cdot,\lambda),v(\cdot,\lambda) = 1$):
\vskip 6pt \noindent\textsc{Lemma 5} ({\it Al-Naggar and Pearson}). {\it Let I $\subset$ $(-\infty,\infty)$ be an interval on which Condition A holds for the general equation (\ref{1.1})}, {\it and let the fundamental system $\{u(\cdot,\lambda),v(\cdot,\lambda)\}$ be defined by the initial conditions (\ref{3.1}) at $x_0 = A$}. {\it Then} \newline
{\bf (i)} {\it There exists a complex valued function $M(\lambda)$ on I which is uniquely defined by the properties:}
\begin{align*}
 & (a) \text{ \ \ \ \ }  Im[M(\lambda)] > 0 \text{ \ \ and \ \ }(b) \text{ \ \ } \lim_{N \rightarrow \infty}  \sfrac{\int_{x_0}^N  (u(x,\lambda)+ M(\lambda) v(x,\lambda))^2  \, dx.}
                        {\int_{x_0}^N  |(u(x,\lambda)+ M(\lambda) v(x,\lambda)|^2 \, dx} = 0.  
\end{align*}
{\bf (ii)} {\it For $\lambda \in I$ the function M in {\bf (i)} is the boundary value of the Titchmarsh-Weyl m-fuction defined by (\ref{4.6}), that is,}
\[ {\it M(\lambda)= \lim_{\epsilon \downarrow 0}} [ m_{A}(\lambda + i\epsilon)]. \]
% {\bf (iii)} For the Sturm-Liouville problem (\ref{4.1})-(\ref{4.2}) the spectral density function is absolutely continuous on $(0,\infty)$ and  
%  \[ \sfrac{1}{\pi} \alpha_{2}(\lambda) =  \lim_{\epsilon\downarrow0}\frac{1}{\pi} Im[m_{A}(\lambda+i\epsilon)]d\mu = f_{A}(\lambda). \] \label{4.9}
%where $\alpha$ is defined in (\ref{4.24}).
\vspace{0.1in}

\noindent\textsc{Proof.} Statements {\bf (i)} and {\bf (ii)} are, respectively, Lemma 1 and Theorem 2 from \cite{AlN1}.\quad\rule{2.2mm}{3.2mm}
\vspace{0.1in}

We now apply these results to the problem (\ref{4.1})-(\ref{4.2}).
\vskip 6pt \noindent\textsc{Lemma 6.}  Assume that for the problem (\ref{4.1})-(\ref{4.2}) {\bf Assumption 3} holds.  Then with $\alpha(\lambda)$ defined by (\ref{4.22})
% the Sturm-Liouville problem (\ref{4.1})-(\ref{4.2}) on $[A,\infty)$ has a.c. spectrum in $(0,\infty)$ and for all $\lambda \in (0,\infty
we have for all $\lambda \in (0,\infty)$
\begin{equation}\label{4.26}
 \alpha(\lambda) =  \lim_{\epsilon \downarrow 0} [ m_{A}(\lambda + i\epsilon) ]. \quad\rule{2.2mm}{3.2mm}  
\end{equation}

\noindent\textsc{Proof.}  By the uniqueness of $M(\lambda)$, and the fact for all $\lambda \in I = (0,\infty)$  $\alpha(\lambda)$ has positive imaginary part (see (4.22) and Lemma 3(ii)) and $\psi_{A}(x,\lambda)$ satisfies  property (b) in Lemma 5(i)  (see Lemma 4)  it follows from Lemma 5(ii) with $I = (0,\infty)$ that for all $\lambda \in (0,\infty)$ we must have (\ref{4.26}).\quad\rule{2.2mm}{3.2mm}
\vspace{0.1in}

We are now ready to prove the representation of $f_{A}(\lambda)$ in the form (\ref{1.5}). 
\vskip 6pt \noindent{\bf Theorem 1} Assume the potential $q$ is given as in (\ref{4.1}) and that {\bf Assumption 3} holds.
Let $\alpha(\lambda)$ be defined as in (\ref{4.22}) and $\psi_{A}(\cdot,\lambda)$ as in (\ref{4.23}). Then the spectral function defined by (\ref{4.8}) for the problem (\ref{4.1})-(\ref{4.2}) is absolutely continuous for $\lambda \in (0,\infty)$ and the corresponding spectral density function admits the following representations for $\lambda \in (0,\infty)$:
\begin{align}
f_{A}(\lambda) & := \rho'(\lambda) = \lim_{\epsilon\downarrow0}\frac{1}{\pi} Im[m(\lambda+i\epsilon)] \label{4.27}\\
& = \sfrac{\alpha_{2}(\lambda)}{\pi} \label{4.28} \\
& = \sfrac{1}{\pi \tilde a(\lambda)}  \label{4.29} \\
& = \sfrac{1}{\pi [P_{1}(x,\lambda) (v(x,\lambda))^{2} + Q_{1}(x,\lambda)v(x,\lambda)v'(x,\lambda) + R_{1}(x,\lambda)(v'(x,\lambda))^{2}]} \label{4.30}
\end{align}
\noindent\textsc{Proof:} The statement (\ref{4.28}), and the absolute continuity of the spectral function $\rho$,
follows from Lemma 6 by taking imaginary parts on each side of equation (\ref{4.26}).  The statement (\ref{4.29}) follows from the definition of $\alpha_2$ in (\ref{4.22}).  To obtain (\ref{4.30}) substitute $U_1 = (P_1,Q_1,R_1)^T$ from the representation (\ref{4.11}) in terms of $\{u(x,\lambda), v(x,\lambda) \}$ and collect coefficients of $\tilde a(\lambda)$, $\tilde b(\lambda)$ and $\tilde c(\lambda)$ to obtain  
\[
P v^{2} + Q v v^{\prime}+ R (v^{\prime})^{2}= \tilde a(\lambda) [W_{x}(v,u)]^{2} = \tilde a(\lambda).
\]
\quad\rule{2.2mm}{3.2mm}

%%%%%Finally, (\ref{4.32}) follows on evaluation of (\ref{4.31}) at $x=A$ using the initial conditions in (\ref{3.1}). \

\vspace{0.1in}
{\bf Remark.} Putting $x=0$ in (\ref{4.30}) yields  $f_{A}(\lambda) = 1/ (\pi R_{1}(A,\lambda))$ which was a main result of Theorem 2 in \cite{AlN2}.
\vspace{0.1in}

The following well known spectral density function formula, due to Titchmarsh \cite{TITCH}, 1946, and Weyl \cite{WEYL}, 1910, in the case of the assumption (\ref{4.3}) (and due to Pearson \cite{PEARSON,FPP1} in the case of assumption (\ref{4.4})) also follows readily from {\bf Theorem 1}: 
\vskip 6pt \noindent{\bf Corollary.} Under the assumptions of {\bf Theorem 1} we have for all $\lambda \in (0,\infty)$ 
\begin{equation} \label{4.31}
        f_A(\lambda) = \lim_{x\to\infty} \frac{1}{\pi [\sqrt{\lambda} (v(x,\lambda)){2} + \frac{1}{\sqrt{\lambda}}(v'(x,\lambda))^{2}]} .
\end{equation}
\vspace{0.1in}

\noindent\textsc{Proof:} The {\bf Assumption 3} ensures (see \cite[Thm 2]{FPP2}) that the solutions $v$ and $u$
 defined by the initial conditions (\ref{3.1}) are bounded for sufficiently large $x$.
 Hence, making use of the initial condition (\ref{1.7}) which $U_1$ satisfies, we have
\begin{align*}
& P_{1}(x,\lambda) (v(x,\lambda))^{2} + Q_{1}(x,\lambda)v(x,\lambda)v'(x,\lambda) + R_{1}(x,\lambda)(v'(x,\lambda))^{2}\\
& = \lim_{x\to\infty} [P_{1}(x,\lambda) (v(x,\lambda))^{2} + Q_{1}(x,\lambda)v(x,\lambda)v'(x,\lambda) + R_{1}(x,\lambda)(v'(x,\lambda))^{2}]\\
& = \lim_{x\to\infty} Q_{1}(x,\lambda)v(x,\lambda)v'(x,\lambda) + \lim_{x\to\infty} [P_{1}(x,\lambda) (v(x,\lambda))^{2} + R_{1}(x,\lambda)(v'(x,\lambda))^{2}]\\
& = 0 + \lim_{x\to\infty} [\sqrt{\lambda} (v(x,\lambda))^{2} + \sfrac{1}{\sqrt{\lambda}} (v'(x,\lambda))^{2}],
\end{align*}
so it follows that (\ref{4.30}) gives rise to the characterization (\ref{4.31}). \quad\rule{2.2mm}{3.2mm}

{\bf Remark.} In \cite{FPP2} we made use of the formula (\ref{4.31}) to establish (\ref{4.30}). Here, by linking the spectral density function first to the m-function, and following the approach of Al-Naggar and Pearson, we have obtained a direct proof of (\ref{4.30}), from which the older result (\ref{4.31}) follows as a consequence.

%%%%%%%%%%%%%%%%%%%%%%%%%%%%%%%%%%%%%%%%%%%%%%%%%%%%%%%%%%%%%%%%%%%%%%%%%%%%%
%%%%%%%%%%%%%%%%%%%%%%%%
%  BELOW: to revise last part of Section 5 -- with refs to Section 4
%%%%%%%%%%%%%%%%%%%%%%%%%%%%%%%%%%%%%%%%%%%%%%%%%%%%%%%%%%%%%%%%%%%%%%%%%%

%This is Section5
\section{Generalization of a Spectral Density Function Characterization to Doubly Singular problems}
\setcounter{equation}{0}
%\begin{equation}
%\rho(\lambda)-\rho(\lambda_{0})=\lim_{\epsilon\rightarrow0}\frac{1}{\pi}
%\int_{\lambda_{0}}^{\lambda}Im[m(\mu+i\epsilon)]d\mu.
%\end{equation}
%
%\begin{equation}
%\label{2_14}\rho(\lambda) = \lim_{\epsilon\searrow0} \int_{0}^{\lambda}
%\mbox{ Im } m(\mu+ i \epsilon) \, d \mu/ \pi
%\end{equation}

In this section we consider the Sturm-Liouville problem
\begin{equation}
\tau y := -y^{\prime\prime}(x)+\left(  \frac{q_{0}}{x^{2}}+\frac{q_{1}}{x^{{}}}%
+\sum_{n=0}^{\infty}\,q_{n+2}x^{n}\text{ }\right)  \text{\ }y(x)=\lambda
y(x),\text{ \ \ \ \ \ }x\in(0,\infty),\label{5.1}%
\end{equation}
\begin{equation}
W_{0}(\left(y(\cdot,\lambda),\phi(\cdot,0)\right) =\lim_{x\rightarrow 0} W_{x}(\left(y(\cdot,\lambda),\phi(\cdot,0)\right)=0, \qquad \text{if $x=0$ is {\bf LC} },\label{5.2}%
\end{equation}
where $\phi(x,0)$ is the first Frobenius solution for $\lambda=0$ given in (\ref{2.16}). Since this is also the principal solution at $x=0$ in all ({\bf LC} and {\bf LP}) cases, the boundary condition (\ref{5.2}) is the Friedrichs boundary condition in the
{\bf LC} cases at $x=0$ and selects $\phi(x,\lambda)$ for all $\lambda \in C$; this boundary condition is also automatically satisfied by $\phi(x,\lambda)$ in all the {\bf LP} cases at $x=0$. In this section we adopt the {\bf Assumption 1} and {\bf 2} from Section 2, and the {\bf Assumption 3} from Section 4.
%in addition the following assumption (see \cite{FPP2,FPP4}):
%\begin{center}
%{\bf Assumption 3:} $\mathbf{Near\ x=\infty:}$
%\end{center}
%For some $x_0>0$ we have either
%\begin{equation}
%q \in L_1(x_0,\infty) \label{5.3}
%\end{equation}
%or
%\begin{equation}
%q' \in L_1(x_0,\infty), \qquad q \in AC_{loc}[x_0,\infty),  \qquad \lim_{x\rightarrow\infty}q(x)=0. \label{5.4}
%\end{equation}
%The assumptions (\ref{5.3}) and (\ref{5.4}) were used in \cite{FPP2,FPP4} to obtain spectral density function characterizations of the type (\ref{1.5}) when the left endpoint is regular. Pearson and Al-Naggar  \cite{AlN1} also obtained a spectral density function characterization of this type when the left endpoint is regular, using a different approach. 
Our aim is to extend the spectral density function characterization in {\bf Theorem 1} (under the above 3 assumptions) to obtain the formula (\ref{1.5}) for the doubly singular problem (\ref{5.1})-(\ref{5.2}).
The {\bf Assumption 1} ensures that $x=0$ is a singular point of type {\bf LP/N} or {\bf LC/N}; the {\bf Assumption 2} ensures that $x=\infty$ is a singular point of type {\bf LP/O-N} with cutoff value $\Lambda = 0$; and the {\bf Assumption 3} ensures that we have a.c. spectrum in $(0,\infty)$. The underlying self-adjoint operator associated with (\ref{5.1})-(\ref{5.2}) has the domain
\begin{eqnarray}   \label{5.3}
  D(A) := \bigg\{ f \in L_2(0,\infty) &|& f(x) \in AC_{loc}^1(0,\infty), \text{ \ }\tau f \in L_2(0,\infty), \nonumber \\
& & lim_{x \to 0} W_x(f(\cdot),\phi(\cdot,0)) = 0 \biggr\}
\end{eqnarray}
in the {\bf LC} cases at $x=0$, and 
\begin{equation} \label{5.4}
D(A):= \biggl\{ f \in L_2(0, \infty) \text{ \ \ } | \text{ \ \ } f\in AC_{\text{loc}}^1 (0, \infty), \text{ \ } \tau f \in L_2(0, \infty) \biggr\}
\end{equation}
in the {\bf LP} cases at $x=0$.
The associated eigenfunction expansion theory which obtains the eigenfunction expansion in the form (\ref{1.2}) for both of the above cases was given in \cite{MN} and \cite{FL}, and explicit formulas for the corresponding Titchmarsh-Weyl m-function and associated scalar spectral function were obtained in \cite{FL} for all cases of the special potential
\begin{equation}
  q(x) = \frac{q_0}{x^2} + \frac{q_1}{x^{}},  \qquad  \text{$q_0$ \ \ and $q_1$ \ \ satisfying  \ \ (\ref{2.2}).} \label{5.5} 
\end{equation}
Here the Titchmarsh-Weyl m-function is defined as in \cite{MN,FL} by
\begin{equation}
\Psi(\cdot,\lambda) := \theta(\cdot,\lambda) - m(\lambda)\phi(\cdot,\lambda) \in L_2(x_0,\infty), \qquad x_0 > 0, \text { \ \ } \text{for all $Im(\lambda) \ne 0$}.\label{5.6}
\end{equation}
where $\phi(\cdot,\lambda)$ and $\theta(\cdot,\lambda)$ are the first and second Frobenius solutions normalized as in (\ref{2.16}).
We now repeat some basic information from \cite{FL} concerning the  Titchmarsh-Weyl m-functions defined by (\ref{5.6}).  In the 
{\bf LC} cases at $x=0$ the m-function defined by (\ref{5.6}) is a Nevanlinna function (that is, a function of class
${\bf N_0}$), and therefore admits an integral representation of the form (\ref{4.7}), 
and the inversion integral for $\rho$ in terms of $m$ is a Titchmarsh-Kodaira formula like (\ref{4.8}). The eigenfunction expansion for the problem (\ref{5.1})-(\ref{5.2}) when $x=0$ is {\bf LC} has the form (\ref{1.2}) where $\phi(\cdot,\lambda)$ is the first Frobenius solution in (\ref{2.16}), and $\rho$ is the spectral function obtained from the above m-function. 

In the {\bf LP} cases at $x=0$, 
the m-function defined by (\ref{5.6}) is a generalized Nevanlinna function of class ${\bf N_{\kappa}}$ for some $\kappa \ge 1$ (see \cite[p. 188]{FL}) and it follows using the theory of Krein and Langer \cite{KL} for these functions (see \cite[Thm 3.5 and Lemma 4.1]{FL}) that they admit the representation
\[
m(z)=(1+z^{2})^{n}\int_{-\infty}^{\infty}\left(
\sfrac{1}{t-z}-\sfrac {t}{1+t^{2}}\right) \,d\sigma(t)+\sum_{j=0}^{m}%
\,\alpha_{j}\,z^{j} 
\]
where $n,m \ge 1$, $\alpha_j \in (-\infty,\infty)$, $\alpha_m \ne 0$ if $m > 0$, and where $\sigma$ is a measure on $\mathbb{R}$ satisfying
\[
 \int_{-\infty}^{\infty}\,\sfrac{d\sigma(t)}{1+t^{2}} <\infty. 
\]
The spectral function for the associated self-adjoint operator A in (\ref{5.4}) is then defined interms of $\sigma$ as 
\begin{equation}
\label{5.7}\rho(\lambda):=\int_{0}^{\lambda}\big(1+s^{2}\big)^{n}\,d\,\sigma(s),\quad \lambda\in(-\infty,\infty).
\end{equation}
While the Titchmarsh-Kodaira formula is well known for {$\bf N_0$} functions, it was only recently established in various cases with two {\bf LP} endpoints and simple spectrum by  Gesztesy and Zinchencko \cite{GZ} and Fulton and Langer \cite{FL}. For the case of equation (\ref{5.1}), where  generalized Nevanlinna functions of class {$\bf N_{\kappa}$} arise, we quote this result, and also relate it to the classical real-variable definition of Levitan and Levinson (see \cite[Thm 4.7 and 4.8]{FL}):
\vskip 6pt \noindent{\bf Theorem 2} ({\it Fulton and Langer}).  Consider equation (\ref{5.1}) with $x=0$ of {\bf LP} type and suppose {\bf Assumptions 1 and 2} hold.  If $\lambda$, $\lambda_0$  are points of $\rho$-measure zero (not discrete eigenvalues of the associated self-adjoint operator A), then the spectral function defined by (\ref{5.7}) has the representation
\begin{equation}
\rho(\lambda) - \rho(\lambda_0) =\lim_{\epsilon\downarrow 0} \frac{1}{\pi} \int_{\lambda_{0}}^{\lambda}Im[m(\mu+i\epsilon)]d\mu \label{5.8}
\end{equation} 
\begin{equation}
\hspace{1.5in} =\lim_{b\to\infty}\sum_{\lambda_{jb}\in (\lambda_0,\lambda)\bigcap\sigma(A_{b})}\frac{1}{\int^{b}_{0}|\phi(x,\lambda_{jb})|^2 dx} . \nonumber
\end{equation}
Note: When (\ref{5.1}) is of {\bf LC} type at $x=0$ we have $n=0$ in (\ref{5.7}) and a standard Nevanlinna representation of $m(\lambda)$, for which the inversion formula is also (\ref{5.8}).
\vspace{0.1in}

Here, $A_b$ is the corresponding truncated self-adjoint operator on $(0,b]$ with any regular boundary condition at $x=b$. It follows from \cite[Thm 4.5]{FL} that for the {\bf LP} cases at $x=0$ the function defined by (\ref{5.7}) 
or (\ref{5.8})  is the spectral function which arises in the eigenfunction expansion (\ref{1.2}). For proofs of convergence results and Parseval relation we refer to \cite{MN,FL}. 
\vspace{0.1in}  

\noindent {\bf Remark.} The above theorem justifies the use of the second formula in  (\ref{5.8}) for the computation of the spectral function $\rho$ as was implemented in the software package {\bf SLEDGE}  when both endpoints are of {\bf LP} type. For Sturm-Liouville problems satisfying {\bf Assumptions 1,2} of Section 2, {\bf SLEDGE} \cite[Equa (1.13)]{TOMS98} does in fact normalize the $\phi$-solution as in (\ref{2.6}), (\ref{2.9}) and (\ref{2.12}).
\vspace{0.1in}

If we make the {\bf Assumption 3} in addition to the assumptions 1 and 2 of Section 2, the spectrum is a.c. on $(0,\infty)$. 
In the case when the left endpoint is regular,
it was shown under {\bf Assumption 3} in \cite[Thm 1 and Cor 2]{FPP2} that the initial value problem (\ref{1.6})-(\ref{1.7}) at $x=\infty$ has a unique solution for each of the cases (\ref{5.3}) and (\ref{5.4}); also, that the spectral density function has the form (\ref{1.5}) (see \cite[Cor 4]{FPP2} and \cite[Thm 1]{FPP4}) for each of these cases. The proof required linking the formula
(\ref{1.5}) to the  well known result of Weyl and Titchmarsh (\ref{4.31}) in the Corollary to Theorem 1. An alternative approach to the proof of (\ref{1.5}) when the left endpoint is regular was given by Al-Naggar and Pearson \cite{AlN1,AlN2} as described in Section 4.  We now follow this method of analysis to generalize the above Theorem 1 to the case when both endpoints are singular and all three assumptions hold. 
Since the spectral function is absolutely continuous in $(0,\infty)$, we may differentiate in (\ref{5.8}) to obtain the spectral density function as 
\begin{equation}
f(\lambda) := \rho'(\lambda) = \lim_{\epsilon\downarrow0}\frac{1}{\pi} Im[m(\lambda+i\epsilon)]. \label{5.9}
\end{equation}  
The eigenfunction expansion associated with the underlying self-adjoint operator $A$ when $x=0$ is {\bf LP} has the form (\ref{1.2}) where $\phi(\cdot,\lambda)$ is the first Frobenius solution in (\ref{2.16}), and $\rho$ is the spectral function defined in (\ref{5.7}) or (\ref{5.8}). In both cases ({\bf LC} and {\bf LP} at $x=0$) the spectral density function is given by (\ref{5.9}) in terms of the m-function which is defined in (\ref{5.6}) above relative to the suitably normalized Frobenius solutions.

We proceed now to transform (\ref{5.9}) to the form (\ref{1.5}); after some analysis this yields (\ref{5.30}) in Theorem 3 below.
Since many of the necessary lemmas which are required are the same as in Section 4, we list and prove only those lemmas which undergo some change as a result of allowing the left endpoint $x=0$ to be a singular endpoint satisfying {\bf Assumptions 1,2}. We also adopt the notational convention that Lemma $n^{*}$ in this section represents the analog of Lemma n in Section 4. 
% where (under {\bf Assumption 3}) $(P,Q,R)^T$ is the unique solution of the initial value problem (\ref{1.6})-(\ref{1.7}). To this end we need the following five lemmas.  It can be shown that the  third component $R(x, \lambda)$ of a solution of Appell's system (\ref{1.6}) satisfies the third order equation of \cite[p. 6584,Equa. 5]{AlN2} (compare also
%\cite[x]{APPELL}):
%\begin{equation}
%\frac{d^{3}R}{dx^{3}} + 4(\lambda-
%q(x)) \frac{\displaystyle dR}{\displaystyle dx} - 2 q^{\prime}(x) R = 0. \label{5.17}
%\end{equation}
%Here, in addition to {\bf Assumption 3} we require that $q \in C^1(0,\infty)$. Also, a general solution of (\ref{5.17}) can be obtained from a suitably normalized solution $R_{1}(x, \lambda)$ as in \cite[p. 6587, Lemma~2]{AlN2}. Here we generalize
%this technique to obtain from a given solution $U_{1}$ of (\ref{1.6}), two other linearly independent solutions. This is the content of the following lemma.

In the sequel it will be helpful to make use of the fundamental system of solutions of Appell's system
(\ref{1.6}) given in (\ref{3.3}) in terms of the suitably normalized Frobenius solutions $\{\phi(\cdot,\lambda), \theta(\cdot,\lambda)\}$ of equation (\ref{5.1}) fixed by the definition (\ref{2.16}). 

\vskip 6pt \noindent\textsc{Lemma $1^{*}$.}  Assume that for equation (\ref{5.1}) {\bf Assumption 1, CaseI,} and {\bf Assumption 2,3} hold.\newline
{\bf (i)} Same as Lemma 1(i); this yields
%For all $\lambda \in (0,\infty)$ let $U_1$ be
%be the unique solution defined at $x=\infty$ in (\ref{4.5}). Then using the indefinite inner product on the solution space of (\ref{1.6}) defined in (\ref{3.7}) we have
\begin{equation}\label{5.10}
\langle U_1, U_1 \rangle= 4 P_1(x, \lambda) R_1(x, \lambda) - (Q_1(x,\lambda))^{2} = 4.
\end{equation}
{\bf (ii)} If the solution $U_1$ is represented in the form (\ref{3.6}), say
\begin{equation}
U_1(x,\lambda)=\left[
\begin{array}
[SLemma 3c]{c}%
P_1(x,\lambda)\\
Q_1(x,\lambda)\\
R_1(x,\lambda)
\end{array}
\right]  =a^{*}\left[
\begin{array}
[c]{c}%
(\theta^{\prime})^{2}\\
-2\theta\theta^{\prime}\\
\theta^{2}%
\end{array}
\right]  +b^{*}\left[
\begin{array}
[c]{c}%
\theta^{\prime}\phi^{\prime}\\
-[\theta\phi^{\prime}+\theta^{\prime}\phi]\\
\theta\phi
\end{array}
\right]  +c^{*}\left[
\begin{array}
[c]{c}%
(\phi^{\prime})^{2}\\
-2\phi\phi^{\prime}\\
\phi^{2}%
\end{array}
\right]  \label{5.11}
\end{equation}
then
\begin{equation}
4 a^{*} c^{*} - (b^{*})^{2} = 4. \label{5.12}
\end{equation}

\noindent\textsc{Proof.} {\bf (ii)} The conversion of the indefinite inner product in terms of the
coefficients $\{a^{*},b^{*},c^{*}\}$ is the property (\ref{3.7}) and (\ref{3.8}); this is readily proved
by substitution of the components of $U_1$ into the inner product formula and use of the wronskian
relation  $W_x(\phi(\cdot,\lambda),\theta(\cdot,\lambda))=1$.
\quad \rule{2.2mm}{3.2mm}

\vskip 6pt \noindent\textsc{Lemma $3^{*}$.}   Assume that for equation (\ref{5.1}) {\bf Assumption 1, CaseI,} and {\bf Assumption 2,3} hold. \newline
{\bf (i)}  Same as Lemma 3(i).\newline
{\bf (ii)} The representation (\ref{5.11}) for $U_1$ has $a^{*}(\lambda) > 0$ and $c^{*}(\lambda) > 0$ for all $\lambda \in (0,\infty)$. 
\vspace{0.1in}

\noindent\textsc{Proof:} {\bf (ii)} In the representation (\ref{5.11}) for $U_1$ we have
\[
R_{1}(x, \lambda) = a^{*} \theta(x, \lambda)^{2} + b^{*} \theta(x, \lambda) \phi(x,
\lambda) + c^{*} \phi(x, \lambda)^{2}.
\]
To see that $a^{*}>0$ and $c^{*}>0$ for all $\lambda\in(0, \infty)$ we first observe
from (\ref{5.12}) that $4a^{*}c^{*} - (b^{*})^{2} = 4$ requires $a^{*}$ and $c^{*}$ to be of the
same sign; and this must hold for all $\lambda \in (0,\infty)$ since they are continuous and cannot pass through zero (which would violate (\ref{5.12})). Factoring out $a^{*}$ we have
\[
R_{1}(x, \lambda) = a^{*} [\theta(x, \lambda)^{2} + (b^{*}/a^{*}) \theta(x, \lambda)
\phi(x, \lambda) + (c^{*}/a^{*}) \phi(x, \lambda)^{2}] 
\]
with the coefficients of $\theta^{2}$ and $\phi^{2}$ positive. Therefore
$R_1$ admits a factorization of the form
\begin{align}
R_{1}(x, \lambda)  & = a^{*}[(\theta(x, \lambda) - \xi\phi(x, \lambda))
(\theta(x, \lambda) - \bar\xi\phi(x, \lambda))]\nonumber\\
& = a^{*} |\theta(x, \lambda) - \xi\phi(x, \lambda)|^{2},\label{5.13}%
\end{align}
where $\xi= \xi_{1} + i \xi_{2}$ must satisfy (by (\ref{5.12}))
\begin{equation}
 \xi_{1} = -b^{*} / (2a^{*}) \mbox{ and } \xi_{2}^{2} = (c^{*} / a^{*}) - [(b^{*})^{2} / (4(a^{*})^{2})] = 1 / (a^{*})^{2}.\label{5.14} 
\end{equation}
Since $a^{*}$, $b^{*}$, and $c^{*}$ are real, we must have either $\xi_{2} = 1/a^{*}$ or
$\xi_{2} = -1 / a^{*}$; but either way the factorization remains the same with
$\xi$ and $\bar\xi$ switched. Since $|\theta-\xi\phi|^{2}>0$, it
follows from the above factorization that $a^{*}(\lambda)>0$ for all $\lambda
\in(0, \infty)$; otherwise, $\lim_{x \to\infty} R_{1}(x, \lambda) \le0$
contradicting the fact that $R_{1}(\infty, \lambda) = 1 / \sqrt{\lambda} > 0$.
Hence, also $c^{*}(\lambda) > 0$ for all $\lambda\in(0, \infty)$. \quad
\rule{2.2mm}{3.2mm}

\vskip 6pt Since the definition of $\xi_{2}$ must be $\pm1 / a^{*}(\lambda)$
from (\ref{5.14}) and since this indeterminacy is actually immaterial, we
choose to take $\xi_{2} = 1 / a^{*}(\lambda)$, so that
\begin{equation}
\xi(\lambda) := -\frac{ b^{*}}{ 2a^{*}} + i \;
\frac{ 1}{a^{*}}. \label{5.15}
\end{equation}
Using the Frobenius solutions normalized by (\ref{2.16}) we now define the complex-valued solution of (\ref{5.1}) for real $\lambda
\in(0, \infty)$, by 
\begin{equation}
\label{5.16}\psi(x, \lambda) := \theta(x, \lambda) - \xi(\lambda) \phi(x,\lambda).
\end{equation}

We now prove the following lemma that equation (\ref{5.1}) satisfies {\em Condition A} for all $\lambda \in (0,\infty)$.
\vskip 6pt \noindent\textsc{Lemma $4^{*}$.} Assume $q(x)$ satisfies {\bf Assumption 1, Case I,} and {\bf Assumptions 2, 3}. Then for $x_{0} > 0$ and all $\lambda\in(0, \infty)$
\begin{equation}
\label{5.17}\lim_{N \to\infty} \frac{\int_{x_{0}}^{N} \psi(x,
\lambda)^{2} \, dx}{\int_{x_{0}}^{N} |\psi(x, \lambda)|^{2} \,
dx } = 0.
\end{equation}
%This means that $\psi(x, \lambda)$ as defined in (\ref{5.23}) satisfies the
%\textit{Condition~A} in (\ref{4.24})
%of Pearson and Al-Naggar \cite{AlN2}
%for all $\lambda\in(0,\infty)$.
% it follows that the
%self-adjoint operators associated with the doubly singular equation (\ref{5.1}) (in both the {\bf LC} and {\bf LP} cases at 0) have
%a.c. spectrum in $(0,\infty)$.

\noindent\textsc{Proof:} Since equation (\ref{4.13}) is satisfied by all
linear combinations of $\phi(\cdot, \lambda)^{2}$, $\theta(\cdot, \lambda)
\phi(\cdot, \lambda)$ and $\theta(\cdot, \lambda)^{2}$ (see \cite[Lemma~1, p.
6584]{AlN2}), it follows that $\psi(\cdot, \lambda)^{2}$ is a solution of
(\ref{4.13}) and since it is complex-valued, also that Re [$\psi(\cdot,
\lambda)^{2}$] and Im [$\psi(\cdot, \lambda)^{2}$] satisfy (\ref{4.13}).
Accordingly, it follows that there exist solutions $U_{2} = (P_{2}, Q_{2},
\mbox{Re}(\psi^{2}))^T$ and $U_{3} = (P_{3}, Q_{3}, \mbox{Im} (\psi^{2}))^T$ of the
first order system (\ref{1.6}), since any real solution $R$ of (\ref{4.13})
can be used to generate a solution of (\ref{1.6}) having $R$ as its third
component; e.g., let $Q^{\prime}= -R^{\prime}$ and $P = - Q^{\prime}/ 2 -
(\lambda- q)R$. From (\ref{5.13}) and (\ref{5.16}) we readily deduce the
following representations of $R_{1}(\cdot, \lambda)$, $Re[ \psi(\cdot,
\lambda)]^{2}$, and $Im[ \psi(\cdot, \lambda)]^{2}$ of the form (\ref{5.11}) (or, the third component of (\ref{5.11})):
\begin{align*}
R_{1}(x, \lambda)  & = a^{*} \theta^{2} - 2a^{*} \xi_{1} \theta\phi+ a^{*} (\xi_{1}^{2} + \xi_{2}^{2}) \phi^{2}\\
\mbox{Re }[ \psi(x, \lambda)]^{2}  & =  \theta^{2} - 2 \xi_{1} \theta\phi+
(\xi_{1}^{2} - \xi_{2}^{2}) \phi^{2}\\
\mbox{Im} [\psi(x, \lambda)]^{2}  & = - 2\xi_{2} \theta\phi+ 2 \xi_{1}\xi_{2} \phi^{2}.\\
\end{align*}
It now follows from these formulas that the above solutions $U_{2}$ and $U_{3}$
associated with $Re[ \psi(\cdot, \lambda)]^{2}$ and $Im[ \psi(\cdot, \lambda
))]^{2}$ are orthogonal to $U_1$ in the sense of the inner product defined in (\ref{3.7}), i.e., we have
\[
\langle U_1, U_{2} \rangle= 2 [a^{*} (\xi_{1}^{2} - \xi_{2}^{2}) + a^{*}
(\xi_{1}^{2} + \xi_{2}^{2})] - 4 a^{*} \xi_{1}^{2} = 0,
\]
and
\[
\langle U_1, U_{3} \rangle= 2 [a^{*} (2 \xi_{1} \xi_{2}) +0] - 4a^{*}
\xi_{1} \xi_{2} = 0.
\]
Hence it follows from Lemma~3(i) that for all $x_{0} > 0$
\[
\lim_{N \to\infty} \frac{ \int_{x_{0}}^{N} Re [\psi(x,
\lambda)]^{2} \, dx}{\int_{x_{0}}^{N} R_{1}(x, \lambda) \, dx} =0
\]
and
\[
\lim_{N \to\infty} \frac{ \int_{x_{0}}^{N} Im [\psi(x,\lambda)]^{2} \, dx}{\int_{x_{0}}^{N} R_{1}(x, \lambda) \, dx} =
0,
\]
from which (\ref{5.17}) follows. \quad\rule{2.2mm}{3.2mm}

\vskip 6pt \noindent\textsc{Lemma $5^{*}$} {\it (Al-Naggar and Pearson).} {\it We assume for equation (\ref{5.1})} {\it that {\bf Assumption 1, Case I,}} {\it and {\bf Assumptions 2,3} hold}. {\it Let I $\subset$ $(-\infty,\infty)$ be an interval on which Condition A holds. Then} \newline
{\bf (i)} {\it There exists a complex valued function $M(\lambda)$ on I which is uniquely defined by the properties:}
\begin{align*}
 (a) \text{ \ \ \ } Im[ M(\lambda)] > 0 \text{ \ \ \ and \ \ \ } (b) \lim_{N \rightarrow \infty}  \sfrac{\int_{x_0}^N  (\theta(x,\lambda)-M(\lambda) \phi(x,\lambda))^2  \, dx.}{\int_{x_0}^N |\theta(x,\lambda) - M(\lambda) \phi(x,\lambda)|^2 \, dx} = 0, \text{ \ \ for \ \ all \ } x_0 > 0.
\end{align*}
where $\{\phi(\cdot,\lambda),\theta(\cdot,\lambda)\}$ are the Frobenius solutions defined in (\ref{2.16}).
\vspace{0.1in}

\noindent\textsc{Proof:} The proof is the same as that for Lemma 1 in \cite{AlN1}.
% except for the sign change which results from the normalization of $\phi$ and $\theta$, $W_{x}(\phi(\cdot,\lambda),\theta(\cdot,\lambda)) = 1$; 
The  proof of this lemma does not depend in any essential way on the choice of the fundamental system of equation (\ref{5.1}). \quad\rule{2.2mm}{3.2mm}

Unfortunately, the corresponding statement {\bf (ii)} from Lemma 5 does not carry over immediately to the
 doubly singular problem (\ref{5.1})-(\ref{5.2}) by borrowing information from \cite{AlN1}; particularly, the proof of
 Theorem 2 in \cite{AlN1} makes use of asymptotic behaviour of solutions which are fixed by initial conditions at a
 regular left endpoint, and therefore do not apply to the Frobenius solutions $\{\phi,\theta\}$. The Titchmarsh-Weyl $m$-function (\ref{5.6}) was first introduced in the papers \cite{GZ,MN,FL}:
 and it wasn't  discovered to be a generalized Nevanlinna function in the {\bf LP} case at $x=0$ until the 2010 paper \cite{FL}. A  direct generalization of Theorem 2 of \cite{AlN1} for cases when the left endpoint is singular remains unknown.
 However, we can recover the analogue of part {\bf (ii)} of Lemma 5 by making appeal to the uniqueness result
 in Lemma $5^{*}$(i) and gleaning information on the boundary behaviour of $m(z)$ from known information on the boundary
 behaviour of $m_{A}(z)$. This is the objective of Lemma $6^{*}$. To this end, it will be helpful to introduce notation for the boundary values of the
``regular'' and ``doubly singular'' Titchmarsh-Weyl functions $m_{A}(\lambda)$ and
$m(\lambda)$ and for the corresponding $\Psi$-functions defined in (\ref{4.6}) and (\ref{5.6}).

{\bf Definition.} Associated with the m-functions, $m_{A}(\lambda)$ and $m(\lambda)$, for the 
problem with regular left endpoint and for the doubly singular problem, respectively, we define 
for all $x \in (0,\infty)$ and all  $\lambda \in (0,\infty)$:
\begin{align}
m_{A}^{+}(\lambda)  & := \lim_{\epsilon\downarrow0} m_{A}(\lambda+ i \epsilon)
),\label{5.18}\\
\Psi_{A}^{+}(x,\lambda)  & := \lim_{\epsilon\downarrow0} \Psi_{A}(\lambda+ i \epsilon) = u(x,\lambda) + m_{A}^{+}(\lambda) v(x,\lambda), 
\label{5.19}\\
m^{+}(\lambda)  & :=  \lim_{\epsilon\downarrow0} m(\lambda+ i \epsilon) ,\label{5.20}\\
\Psi^{+}(x,\lambda)  & := \lim_{\epsilon\downarrow0} \Psi(\lambda+ i \epsilon) = \theta(x,\lambda) - m^{+}(\lambda) \phi(x,\lambda),\label{5.21}
\end{align}
%Since for $\epsilon > 0$, we have $Im[m_{A}(\lambda + i\epsilon)] > 0$ and $Im[m(\lambda + i\epsilon)] > 0$ it follows that
%$Im[m_{A}^{+}(\lambda)] > 0$ and $Im[m^{+}(\lambda)] > 0$. 
%The following lemma gives information needed for the ``doubly singular''
%analogue of equation (\ref{4.29}) or Lemma 4(ii) (that is, the analogue of \cite[Thm2]{AlN1}).
\vskip 6pt \noindent\textsc{Lemma $6^{*}$.}  We assume for equation (\ref{5.1}) that {\bf Assumption 1, Case I,} and {\bf Assumptions 2,3} hold. Then \newline
{\bf (i)} for all $x \in(0, \infty)$ and all $z$       with Im $z \ne 0$
\begin{equation}
\label{5.22}\Psi(x, \lambda) = \frac{\Psi_{A}(x,z)}{- \phi^{\prime}(A,z) + m_{A}(z)\phi(A,z)}.
\end{equation}
{\bf (ii)} For all $x \in(0, \infty)$ and all $\lambda\in(0, \infty)$
\begin{equation}
\label{5.23}\Psi^{+}(x, \lambda) = \frac{\Psi_{A}^{+}(x,\lambda)}{- \phi^{\prime}(A, \lambda) + m_{A}^{+}(\lambda)\phi(A,\lambda)}.
\end{equation}
{\bf (iii)} For $x_{0} > 0$ and all $\lambda \in (0,\infty)$
\begin{equation}
\label{5.24}\lim_{N \to\infty} \frac{\int_{x_{0}}^{N} \Psi^{+}(x, \lambda)^{2} \, dx}{ \int_{x_{0}}^{N} |\Psi^{+}(x,\lambda)|^{2} \, dx} = 0.
\end{equation}
%(iv) For all $\lambda\in(0, \infty)$
%\begin{equation}
%\label{3_35}m_{+}(\lambda) = \frac{\displaystyle M_{+}(\lambda) \theta(x_{0},
%\lambda) - \theta^{\prime}(x_{0},\lambda)}{\displaystyle M_{+}(\lambda)
%\phi(x_{0}, \lambda) - \phi^{\prime}(x_{0}, \lambda)}.
%\end{equation}
{\bf (iv)} For all $\lambda\in(0, \infty)$
\begin{equation}
\label{5.25}\mbox{Im } m_{+}(\lambda) = \frac{\displaystyle \mbox{Im }
m_{A}^{+}(\lambda)}{\displaystyle |m_{A}^{+}(\lambda) \phi(A, \lambda)-
\phi^{\prime}(A, \lambda)|^{2}} > 0.
\end{equation}
{\bf (v)} For all $\lambda\in(0, \infty)$ the complex valued function $\xi(\lambda)$ defined in (\ref{5.15}) is the boundary value of the Titchmarsh-Weyl m-function, that is,
\begin{equation} \label{5.26}
\xi(\lambda) = \lim_{\epsilon\downarrow0} m(\lambda+i\epsilon).
\end{equation}.
\vspace{0.1in}

\noindent\textsc{Proof.} For (i) make use of the fact that for $Im z \ne 0$ both $ \Psi_{A}(x,z) = u(x, \lambda) + m_{A}(\lambda) v(x, \lambda)$ and $\Psi(x,z)=\theta(x, \lambda)$ - $m(\lambda) \phi(x,\lambda)$ are in $L_{2}(x_{0}, \infty)$, and therefore linearly dependent. Using
$W_{x}(\phi, \theta) = W_{x}(\phi, \Psi) = 1$, the relation of linear
dependence is found to be
\[
\Psi(x,z) = \frac{\Psi_{A}(x,z)}{W_{x}(\phi, u + m_{A}(v))},
\]
and (\ref{5.22}) follows on evaluation of the denominator at $x=A$ making
use of the initial conditions (\ref{3.1}). For (ii) put $z=\lambda+ i \epsilon$ in (\ref{5.22}) and pass $\epsilon\to0$. For (iii) put
(\ref{5.23}) into (\ref{5.24}) and factor the constant terms out of the
integrals to get the equivalent statement
\[
\lim_{N \to\infty} \frac{\int_{x_{0}}^{N} \Psi_A^{+}(x,\lambda)^{2} \, dx}{\int_{x_{0}}^{N} |\Psi_A^{+}(x, \lambda)|^{2} \,
dx } = 0.
\]
But $\Psi_{A}^{+}(x,\lambda)$ is known to  satisfy {\it Condition A}  by Lemma 5(ii), Lemma 6, and Lemma 4 with (i.e. by Theorem 2 of \cite{AlN1} applied to the problem (\ref{4.1})-(\ref{4.2})). Hence it follows that (\ref{5.24}) also holds for all $\lambda \in (0,\infty)$ and all $x_0 > 0$; in other words, $\Psi^{+}(x,\lambda)$
% (which involves the Frobenius solutions $\theta$ and $\phi$)
also satisfies {\it Condition A}.  For (iv)  use 
\[ m^{+}(\lambda) = W_{x}(\theta(\cdot,\lambda),\Psi^{+}(\cdot,\lambda)) \]
and substitute the right hand side of (\ref{5.23}) for $\Psi^{+}(x,\lambda)$, evaluating  $W_{x}(\theta,u)$ and $W_{x}(\phi,v)$ at $x=A$ using the initial conditions (\ref{3.1}), to obtain for all $\lambda \in (0,\infty)$
\[ m^{+}(\lambda) = \sfrac{\theta(A,\lambda) m_{A}^{+}(\lambda) - \theta^{\prime}(A,\lambda)}{\phi(A,\lambda) m_{A}^{+}(\lambda) - \phi^{\prime}(A,\lambda)} .\]
Then (\ref{5.25}) follows by taking imaginary parts on both sides. The right hand side of (\ref{5.25})
 is positive for all $\lambda \in (0,\infty)$ because the numerator is positive by Theorem 1 (equation (\ref{4.28})), and the denominator never becomes zero (by separating real and imaginary parts and observing that the cases $\phi(A,\lambda)$ = 0 and $\phi(A,\lambda) \ne  0$ both yield a positive denominator).
To prove {\bf (v)} we observe first that for all $\lambda \in (0,\infty)$, $m^{+}(\lambda)$ has positive imaginary part (by (\ref{5.25})) and $\Psi^{+}(x,\lambda)$  satisfies (by (\ref{5.24})) property (b) in Lemma $5^{*}$(i).  Similarly, the function $\xi$ defined in (\ref{5.15}) has positive imaginary part (see Lemma 3(ii) ) for all $\lambda \in (0,\infty)$ and  the solution $\psi(x,\lambda)$ defined in (\ref{5.16}) satisfies (see Lemma $4^{*}$) property (b) in Lemma $5^{*}$(i). Hence by the uniqueness of the function $M(\lambda)$ satisfying the two properties of Lemma $5^{*}$(i), the functions $m^{+}(\lambda)$ and $\xi(\lambda)$ must be identical, that is, (\ref{5.26}) holds. \quad\rule{2.2mm}{3.2mm} 

We are now ready to prove the representation of $f(\lambda)$ in the form (\ref{1.5}). This represents the  ``doubly singular'' analogue of the spectral density function characterization of {\bf Theorem 1}.

\vskip 6pt \noindent{\bf Theorem 3} We assume that {\bf Assumption 1, Case I,} and {\bf Assumptions 2,3} hold.
Let $\xi(\lambda)$ be defined as in (\ref{5.15}) and $\psi(\cdot,\lambda)$ as in (\ref{5.16}). Then the spectral function defined by (\ref{5.8}) for the problem (\ref{5.1})-(\ref{5.2}) (in both the {\bf LC} and {\bf LP} cases at $x=0$) is absolutely continuous for $\lambda \in (0,\infty)$ and the corresponding spectral density function admits the following representations for $\lambda \in (0,\infty)$:
\begin{align}
f (\lambda) & := \rho'(\lambda) = \lim_{\epsilon\downarrow0}\frac{1}{\pi} Im[m(\lambda+i\epsilon)] \label{5.27}\\
& = \sfrac{\xi_{2}(\lambda)}{\pi} \label{5.28} \\
& = \sfrac{1}{\pi a^{*}(\lambda)}  \label{5.29} \\
& = \sfrac{1}{\pi [P_{1}(x,\lambda) (\phi(x,\lambda))^{2} + Q_{1}(x,\lambda)\phi(x,\lambda)\phi'(x,\lambda) + R_{1}(x,\lambda)(\phi'(x,\lambda))^{2}]}. \label{5.30}
\end{align}

\noindent\textsc{Proof:} The statement (\ref{5.28}), and the absolute continuity of the spectral function $\rho$, follows from Lemma $6^{*}$(v) by taking imaginary parts of each side of (\ref{5.26}).  The statment (\ref{5.29}) follows from the definition of $\xi_2$ in (\ref{5.15}).  To obtain (\ref{5.30}) substitute $U_1 = (P_1,Q_1,R_1)^T$ from the representation (\ref{5.11}) in terms of $\{\theta(x,\lambda), \phi(x,\lambda) \}$ and collect coefficients of $a^{*}(\lambda)$, $b^{*}(\lambda)$ and $c^{*}(\lambda)$ to obtain  
\[P (\phi)^{2} + Q \phi \phi^{\prime}+ R (\phi)^{\prime})^{2}=  a^{*}(\lambda) [W_{x}(\phi,\theta)]^{2} =  a^{*}(\lambda). \quad\rule{2.2mm}{3.2mm} \]

%Using (\ref{5.31}) we now have the following ``doubly singular'' analogue of the classical spectral density function formula of Titchmarsh (1946) and Weyl (1910): 
%\vskip 6pt \noindent{\bf Corollary} Under the assumptions of {\bf Theorem 3} we have 
%\begin{equation} \label{5.29}
%f_(\lambda) = \lim_{x \rightarrow \infty}\frac{1}{\pi [\sqrt{\lambda} (\phi(x,\lambda)^{2} + \frac{1}{\sqrt{\lambda}}(\phi'(x,\lambda))^{2}]}  
%\end{equation}
%\vspace{0.1in}
%Ques:  Do we have boundedness of \theta and \phi solutions near \infty ??  (Needed this in proof of previous Corollary)
%\noindent\textsc{Proof.} The proof is the same as the proof of the Corollary of {\bf Theorem 1}. The only new aspect is that the solution
%$\phi(\cdot,\lambda)$ is the suitably normalized Frobenius solution in (\ref{2.16}). \quad\rule{2.2mm}{3.2mm}

\section{Some Examples with Explicitly Known Spectral Density Functions}
\setcounter{equation}{0}
In this section we give the explicit formulas from \cite{MN,FL} for the Frobenius solution $\phi(\cdot,\lambda)$, the Titchmarch-Weyl m-function, and the spectral density function for some examples of the special potential (\ref{5.5}). We restrict attention to those cases which will be used as test problems for the numerical algorithms in Sections 8 and 9.
\vspace{0.1in}

\noindent{\bf Example 1:} [$q_0 = -a, a > 0; q_1 = \ell(\ell+1)$] Hydrogen Atom
\begin{equation}  \label{6.1} 
  -y''+\left( -\frac{a}{x} +\frac{\ell (\ell+1)}{x^2} \right) y = \lambda y 
  \qquad  a > 0, \qquad 0 < x < \infty.
\end{equation}
The first Frobenius solution with normalization (\ref{2.9}) is
\begin{eqnarray}   \label{6.2}
  \phi(x, \lambda) &=& x^{\ell+1}[1 + \sum^{\infty}_{n=1}a_n(\lambda) x^n]
                       \nonumber \\ 
  &=& x^{\ell + 1} e^{ix\sqrt{\lambda}} M(\ell + 1 - \beta, \, 2\ell + 2, \, 
      -2ix\sqrt{\lambda})  \nonumber \\
  &=& \sfrac{1}{(-2i\sqrt{\lambda})^{\ell + 1}} {\cal M}_{\beta, \ell + 
      \frac{1}{2}}(-2ix\sqrt{\lambda}),
\end{eqnarray}
with $\beta := ia / 2 \sqrt{\lambda}$ for all $\lambda \in \Complex$. Here M is the confluent
hypergeometric function of first kind and $\cal M$ is the corresponding Whittaker function of first kind.
The coefficients $a_n(\lambda)$ are polynomials in $\lambda$ of degree $[n/2]$ which are generated from
the recurrence relation 
\[
a_n(\lambda)=-\frac{a}{n(n+2\ell+1)}a_{n-1}(\lambda)-\frac{\lambda}{n(n+2\ell+1)}a_{n-2}(\lambda),
\]
and the first three are
\begin{eqnarray*}
a_1&=&-\frac{a}{2\ell+2},\nonumber\\
a_2(\lambda)&=&\frac{a^2-2(\ell+1)\lambda}{2!(2\ell+2)(2\ell+3)},\nonumber\\
a_3(\lambda)&=&\frac{-a^3+(6\ell+8)a\lambda}{3!(2\ell+2)(2\ell+3)(2\ell+4)},\nonumber
\end{eqnarray*}
The Titchmarsh-Weyl m-function arising from (\ref{5.6}) is
\begin{equation}   \label{6.3} 
   m_\ell (\lambda ) = k_{\ell}(\lambda) \left[-a \log(-2i\sqrt{\lambda})
   -a \Psi(1 - ia / (2 \sqrt{\lambda})) - 2 \gamma a + i \sqrt{\lambda} 
    \right] + p_\ell(\lambda),
\end{equation} 
where $\Psi$ is the psi or digamma function, $\gamma$ is Euler's
constant, 
\[  k_{\ell}(\lambda) := \sfrac{1}{[(2\ell + 1)!]^2} \prod^{\ell}_{j=1}(4 \lambda j^2 + a^2), \]
and where $p_{\ell}(\lambda)$ is a polynomial of degree $\ell$ (see \cite{MN,KF}). We take $0 \le arg(\lambda) < 2\pi$, so that the branch cut for $\sqrt{\lambda}$ and $m_{\ell}$ is on the positive real $\lambda$-axis. \newline

\noindent The associated spectral density function arising from (\ref{5.9}) is
\begin{equation} \label{6.4}
   f_\ell (\lambda) := \lim_{\epsilon \to 0} \frac{Im [m_{\ell} (\lambda 
      + i \epsilon)]}{\pi} = k_{\ell}(\lambda) \left[\frac{a}{1-e^{-\pi a / 
      \sqrt{\lambda}}}\right].
\end{equation}
\vspace{0.1in}

\noindent{\bf Example 2:} [$q_0 = -a, a < 0; q_1 = \ell(\ell+1)$]  Repulsive Coulomb
\begin{equation}
\label{6.5}-y^{\prime\prime}+ \left( -\frac{a}{x}
+ \frac{ \ell(\ell+1)}{ x^{2}} \right)  y = \lambda y, \qquad a < 0, \qquad  0 < x , \infty
\end{equation}
The first Frobenius solution with normalization (\ref{2.9}) is (same as (\ref{6.2}) with $a < 0$)
\begin{eqnarray}
  \phi(x, \lambda) &=& x^{\ell+1}[1 + \sum^{\infty}_{n=1}a_n(\lambda) x^i]
                       \nonumber \\ 
  &=& x^{\ell + 1} e^{ix\sqrt{\lambda}} M(\ell + 1 - \beta, \, 2\ell + 2, \, 
      -2ix\sqrt{\lambda})  \nonumber \\
  &=& \sfrac{1}{(-2i\sqrt{\lambda})^{\ell + 1}} {\cal M}_{\beta, \ell + 
      \frac{1}{2}}(-2ix\sqrt{\lambda}), \label{6.6}
\end{eqnarray}
where $\beta := ia / 2 \sqrt{\lambda}$, and $a < 0$. The recurrence relation and first three coefficients, $a_n(\lambda)$, are the same as in Example 1 with $a < 0$.
The Titchmarsh-Weyl m-function arising from (\ref{5.6}) is
\begin{equation}   \label{6.7} 
   m_\ell (\lambda ) = k_{\ell}(\lambda) \left[-a \log(-2i\sqrt{\lambda})
   -a \Psi(1 - ia / (2 \sqrt{\lambda})) - 2 \gamma a + i \sqrt{\lambda} 
    \right] + p_\ell(\lambda).
\end{equation} 
where $a < 0$, and $k_{\ell}(\lambda)$,  $p_{\ell}(\lambda)$ are the same, with $a < 0$, as given
above for the hydrogen atom. The branch cut is taken again on the positive real $\lambda$-axis. \newline

\noindent The associated spectral density function arising from (\ref{5.9}) is
\begin{equation}
\label{6.8}f_{\ell}(\lambda) = k_{\ell}(\lambda)\left( \frac{ |a|}{ \exp(|a| \pi/ \sqrt{\lambda}) - 1} \right),
\end{equation}
\vspace{0.1in}

\noindent{\bf Example 3:} [$q_0 = \nu^2 - 1/4, \nu \ne M/2, M = 0,1,2 \cdots; q_1 = 0 $]  Bessel Equation of Non-integer Order
\begin{equation}
\label{6.9}-y^{\prime\prime}+ \frac{ \nu^{2} - 0.25}{ x^{2}} \, y = \lambda y, \qquad\nu> 0, \; \nu\ne N / 2, N = 1,2, \ldots,
\end{equation}
The first Frobenius solution with normalization (\ref{2.6}) is
\begin{equation}
\label{6.10}  \phi(x,\lambda):=  x^{\nu+0.5} \left[ 1 + \sum_{j=1}^{\infty} \sfrac{(-1)^j \lambda^j x^{2j}}{j! (\nu+1)_j 2^{2j}} \right] =  2^{\nu} \Gamma(\nu+1) \lambda^{-\nu/2}x^{1/2} J_{\nu}(\sqrt{\lambda} x)
\end{equation}
The Titchmarsh-Weyl m-function arising from (\ref{5.6}) is
\begin{equation}\label{6.11}
 m(\lambda) = -\sfrac{\pi}{2^{2\nu+1} \Gamma^2(\nu+1)\cdot sin(\nu \pi)} e^{-i\nu\pi} \lambda^{\nu} 
\end{equation}
where $0 \le arg(\lambda) < 2\pi$, so that the branch cut for $\lambda^{\nu}$ is on the positive real $\lambda$-axis.\newline

\noindent The associated spectral density function arising from (\ref{5.9}) is
\begin{equation}
\label{6.12}f_{\nu}(\lambda) = \frac{\lambda^{\nu}}{ 2^{2 \nu+ 1} \Gamma^{2}(\nu+ 1)}.
\end{equation}
\vspace{0.1in}

\noindent{\bf Example 4:} [$q_0 = N^2 - 1/4, N = 0,1,2, \cdots; q_1 = 0$]  Bessel Equation of Integer Order
\begin{equation}
\label{6.13}-y^{\prime\prime}+ \frac{\displaystyle N^{2} - 0.25}
{\displaystyle x^{2}} \, y = \lambda y, \qquad a > 0, N = 0, 1, \ldots,
\end{equation}
 The first Frobenius solution with normalization (\ref{2.12}) is
\begin{equation}
\label{6.14}  \phi(x,\lambda):=  x^{N+0.5} \left[ 1 + \sum_{j=1}^{\infty} \sfrac{(-1)^j \lambda^j x^{2j}}{j! (N+1)_j 2^{2j}} \right] =  2^{\nu} \Gamma(N+1) \lambda^{-N/2}x^{1/2} J_{N}(\sqrt{\lambda} x), N = 0,1,2,\cdots .
\end{equation}
The Titchmarsh-Weyl m-function arising from (\ref{5.6}) is
\begin{align}
m_{0}(\lambda)  & = -\log(-2i\sqrt{\lambda}) + \gamma -2\ell n 2 \nonumber \\
m_{N}(\lambda)  & = \frac{\lambda^{N}}{2^{2N}
(N!)^{2}} m_{0}(\lambda) + \frac{\lambda^{N} H_{N+1}}{ 2^{2N+1}(N!)^{2}}, \qquad N \ge 1. \label{6.15}
\end{align}
where $0 \le arg(\lambda) < 2\pi$, so that the branch cut for $m_{0}$ is on the positive real $\lambda$-axis.\newline

\noindent The associated spectral density function arising from (\ref{5.9}) is
\begin{equation}
\label{6.16} f_N(\lambda) = \frac{\lambda^N}{2^{2N+ 1}(N!)^2}.
\end{equation}

\section{Numerical Methods}                                           % 5
\setcounter{equation}{0}

In this section and the following two sections we describe some new numerical methods for obtaining 
approximations to the spectral density function (\ref{1.4}) and then the spectral function, by making use of the
new representation (\ref{5.30}) in Theorem 3, and compare their performance with {\bf SLEDGE}. For general
information and discussion of numerical methods for Sturm-Liouville problems we refer to Pryce's book \cite{PRYCE}
and for spectral function computation using {\bf SLEDGE} we refer to our previous papers \cite{SLEDGE,EAST,TOMS98}. 

Many numerical methods for (\ref{1.1}) break down near a singular point
at $x=0$.  However, when we take this singular point to be a regular singular point, it admits a
convergent Frobenius expansion, and then a finite number of leading terms 
in the sum can be used as an initial approximation near $x=0$. For equation (\ref{2.4}) the 
indicial equation is (\ref{2.5}) and the principal solution is the first Frobenius solution
with the larger indicial root,
\[   r_1 = \nu: = 0.5 + 0.5 \sqrt{1 + 4 q_{0}}.   \]
This solution has the general form
\begin{equation}   \label{7.1}
   \phi(x,\lambda) = \sum_{n=0}^{\infty} a_n x^{n + \nu},
\end{equation}
and the general recurrence formula is  
\[  a_1 = \sfrac{q_{1} a_0}{\nu (\nu + 1) - q_{0}},  \]
and for $n > 1$
\[
  a_n = \sfrac{ - \lambda a_{n-2} + q_1 a_{n-1} +  \sum_{k=0}^{n-2} q_{k+2} a_{n-2-k}}
          {(\nu + n -1)(\nu + n) - q_{0}}.  
\]
The choice for $a_0$ fixes the normalization of $\phi(\cdot,\lambda)$, and in this paper 
we have made the simple choice $a_0 = 1$ in all the cases (\ref{2.6}), (\ref{2.9}) and (\ref{2.12}); this
ensures that the properties {\bf (i),(ii),(iii)} of section 2 hold in all the cases of {\bf Assumption 1}.  

We note that there is a risk of loss of significance in the computation of $a_n(\lambda)$ for very large 
$\lambda$, because for moderate $n$ the powers of $\lambda$ in the numerator of $a_n$ build up 
faster than the denominator does.  We have found that keeping
\vspace{0.1in}
\begin{equation} \label{7.2}
x < x_0 (\lambda) : = |q_{0}| / \sqrt{\lambda} 
\end{equation}
and using the truncated Frobenius series (to full machine precision) only on $(0, x_0(\lambda)]$ 
works well.  For $x > x_0 (\lambda)$ we use the methods from \cite{FPP2} 
and \cite{FPP4} for regular problems.  In brief, the algorithm is as follows:

\noindent 
{\bf(i)} For a given $\lambda$ choose a `matching point' $x(\lambda)$.

\noindent 
{\bf (ii)} Use the first Frobenius solution (\ref{7.1}) (this is the principal solution near zero) on the
interval $(0, x_0 (\lambda)]$ to produce values for solution and its derivative (usually to machine precision) at $x_0 (\lambda)$.

\noindent
{\bf (iii)} Apply standard methods to numerically estimate $y$ for (\ref{5.1}) 
on the interval $[x_0(\lambda), x(\lambda)]$ satisfying initial conditions from {\bf (i)}.  
Because of the oscillation in $y$ for $\lambda > 0$, we use  a piecewise 
trigonometric approximation to $y$.
% over a piecewise polynomial one.

\noindent
{\bf (iv)} Approximate the solution  $P(x, \lambda)$, $Q(x, \lambda)$, and 
$R(x, \lambda)$ of the initial value problem (\ref{1.6})-(\ref{1.7}) at the matching point $x(\lambda)$ 
using one of the approaches described below.

\noindent
{\bf (v)} Substitute the estimates from {\bf (iii)},{\bf (iv)} into 
\begin{equation}  \label{7.3}
  f_x(\lambda) := \sfrac{1}{\pi [P(x, \lambda) y(x, \lambda)^2 + 
    Q(x, \lambda) y(x, \lambda) y'(x, \lambda) + R(x, \lambda) 
    y'(x, \lambda)^2]},
\end{equation}
to produce an approximation to the exact spectral density, $f(\lambda)$, given in Theorem 3, equation (\ref{5.30}).

\vskip 6pt
The two papers \cite{FPP2} and \cite{FPP4} derive two very different
approaches to the numerical computation of $(P(b, \lambda)$, 
$Q(b, \lambda)$, and $R(b, \lambda))$.   In \cite{FPP2} we constructed 
a family of recurrence formulas that generated successively more accurate 
approximations to $P$, $Q$, and $R$ and hence to $f(\lambda)$.  For 
fixed $x > 0$ define for each positive integer $j$ the family of functions
\begin{equation}    \label{PQR}
  F_x^j(\lambda) := \sfrac{1}{\pi \left[P_j y^2 + Q_j y y' + R_j y'^2 \right]}.
\end{equation}
From \cite{FPP1} for $j=1$ we define 
\begin{equation}      \label{PQR1}
   P_1 := \sqrt{\lambda}, \; Q_1 := 0, \; R_1 := 1 / \sqrt{\lambda}.
\end{equation}
From \cite{FPP2} the next formula in the family for $j=2$ is defined by 
\begin{equation}      \label{PQR2}
   P_2 := \sqrt{\lambda - q(x)}, \quad Q_2 := -q'(x) / [2 (\lambda - q(x))^{3/2}],
   \quad R_2 := 1 / \sqrt{\lambda - q(x)}.
\end{equation}
The final member that we use is $j=3$, also defined in \cite{FPP2} as:
% it is more complicated.
\begin{eqnarray}
   P_3 &:=& P_2 + 0.25 \gamma_2 + 0.125 \gamma_1^2 / \gamma_0 \nonumber \\
   Q_3 &:=& Q_2 - \sfrac{d}{dx} \left\{ -0.25 \gamma_0^2 \gamma_2 + 0.125
           \gamma_0 \gamma_2^2  \right\}                  \label{PQR3} \\
   R_3 &:=& R_2 - 0.25 \gamma_2 + 0.125 \gamma_0 \gamma_1^2
           \nonumber
\end{eqnarray}
where
\[    \gamma_k := \sfrac{d^k}{dx^k} \left[ \sfrac{1}{\sqrt{\lambda - q(x)}}
                  \right]  \]
for $k = 0, 1, 2.$  For regular problems on $[A,\infty)$,  $A > 0$, it is 
shown in \cite{FPP2} that each member of this family converges to the 
spectral function $f_{A}(\lambda)$ as $x \to \infty$.   For the hydrogen atom 
potential (\ref{6.1}) on $[A,\infty)$ the theory of \cite{FPP2}  implies that 
\begin{equation}    \label{FPP2rate}
   f_{A}(\lambda) - F_x^j(\lambda) = O(1 / x^{2j-1}) \mbox{ as } x \to \infty.
\end{equation}
This method requires knowledge of derivatives of the potential $q(x)$.

In \cite{FPP4} we constructed explicit approximations to the solutions of 
(\ref{1.6}) with known residual terms that arise during the construction; 
in particular, replace $P(x, \lambda)$, $Q(x,\lambda)$, and $R(x,\lambda)$ with estimates of the 
form
\begin{eqnarray}
  P_N(x) &:=& \sqrt{\lambda} + \sum_{j=1}^N a_j / x^j     \nonumber \\
  Q_N(x) &:=& \qquad \; \; \sum_{j=1}^N b_j / x^{j+1} \label{EXPANSION} \\
  R_N(x) &:=& 1 / \sqrt{\lambda} + \sum_{j=1}^N c_j / x^j,     \nonumber
\end{eqnarray}
where $\{ a_j \}$,  $\{ b_j \}$, and $\{ c_j \}$ will depend on $\lambda$ 
but not $x$.  The resulting sums are substituted into (\ref{1.6}) and the 
coefficients chosen to match terms of like powers.  Then we put the computed 
coefficients into (\ref{EXPANSION}) and define the family of approximations
\begin{equation}   \label{fPQR}
  f_x^N(\lambda) := \frac{1}{\pi\lbrack P_{N}(x,\lambda
)y(x,\lambda)^{2}+Q_{N}(x,\lambda) y(x,\lambda) y^{\prime}(x,\lambda
)+R_{N}(x,\lambda) y^{\prime}(x,\lambda)^{2}]}.
%\sfrac{1}{\pi [(P_N y^2)(x, \lambda) + (Q_N y)(x, \lambda)
%               y'(x, \lambda) + (R_N y'^2)(x, \lambda)]}.
\end{equation}
Specifically, the $N$th residuals are defined as
\vskip 6pt 
\[
\left(
\begin{array}
[c]{l}
\phi_N^P \\
\phi_N^Q\\
\phi_N^R 
\end{array}
\right) :=\left(
\begin{array}
[c]{c}%
P_N^{\prime}\\
Q_N^{\prime}\\
R_N^{\prime}
\end{array}
\right) - \left(
\begin{array}
[c]{ccc}%
0 & \lambda-q & 0\\
-2 & 0 & 2(\lambda-q)\\
0 & -1 & 0
\end{array}
\right)  \cdot\left(
\begin{array}
[c]{c}%
P_N\\
Q_N\\
R_N
\end{array}
\right) 
\]
\vskip 6pt 
\noindent
and we attempt to make these small, as $x \to \infty$, by the choice of 
coefficients in (\ref{EXPANSION}).   All potentials in the examples of 
the previous section have the form (\ref{5.5}), that is, 
\begin{equation}    \label{GENq}
  q(x) = A / x + B / x^2,
\end{equation}
where $A = q_1$ and $B=q_0$ and (\ref{2.2}) is satisfied. It is straightforward to show that
\begin{eqnarray*}
  \phi_N^P &=& \left[ \sum_{j=1}^N \sfrac{-ja_j - \lambda b_j + A b_{j-1} + 
          B b_{j-2}}{x^{j+1}}  \right] + \sfrac{B b_{N-1} + Ab_N}{x^{N+2}} 
        + \sfrac{B b_N}{x^{N+3}}   \\
  \phi_N^Q &=& \left[ \sum_{j=1}^N \sfrac{-(j-1)b_{j-2} + 2a_j - 2 
          \lambda c_j + 2Ac_{j-1} + 2 Bc_{j-2}}{x^j}  
          + \sfrac{2A + 2 B / x}{x\sqrt{\lambda}}  \right]    \\
     & & - \sfrac{N b_{N-1} + 2 A c_N + 2 B c_{N-1}}{x^{N+1}}
         + \sfrac{2 B c_N - (N+1) b_N}{x^{N+2}}  \\
  \phi_N^R &=& \sum_{j=1}^N \sfrac{-j c_j + b_j}{x^{j+1}}.
\end{eqnarray*}
If we require the coefficients to satisfy
\begin{eqnarray}    
 j a_j + \lambda b_j &=& Ab_{j-1} + B b_{j-2}  \label{COND_P} \\
 a_j - \lambda c_j &=& (j-1) b_{j-2} / 2 - Ac_{j-1} - Bc_{j-2}
    + [A \delta_{j1} + B \delta_{j2}] / \sqrt{\lambda}  
     \label{COND_Q} \\
   b_j -j c_j &=& 0,        \label{COND_R}
\end{eqnarray}
for $j = 1, 2, \ldots, N$, then the residuals simplify to
\begin{eqnarray}
  \phi_N^P &=& \sfrac{A b_N + B b_{N-1}}{x^{N+2}} + \sfrac{B b_N}{x^{N+3}} \\
  \phi_N^Q &=& \sfrac{-N b_{N-1} + 2 A c_N + 2 B c_{N-1}}{x^{N+1}}
               +\sfrac{2 B c_N - (N+1) b_N}{x^{N+2}}     \\
  \phi_N^R &=& 0.
\end{eqnarray}
If we adopt the convention that coefficients with nonpositive subscripts 
have zero values, then the solution of (\ref{COND_P})--(\ref{COND_R}) can 
be written, for $1 \le j \le N$,
\begin{eqnarray*}
  a_j &=& (t_1 + t_2) / 2  \\
  c_j &=& (t_1 - t_2) / (2 \lambda) \\
  b_j &=& j c_j,
\end{eqnarray*}
where
\[ t_1 = [(A b_{j-1} + B b_{j-2}] / j  \]
and
\[ t_2 = 0.5 (j-1) b_{j-2} - A c_{j-1} - B c_{j-2}
    - [A \delta_{j1} + B \delta_{j2}] / \sqrt{\lambda}. \]
Since the derivatives of the residuals do not change sign once $x$ is 
sufficiently large, the theory of \cite{FPP2} implies for that
\begin{equation}      \label{FPP4rate}
  f(\lambda) - f_x^N(\lambda) = O(1 / x^{N+1})
\end{equation}
as $x \to \infty$.

To numerically estimate the spectral density function 
$f(\lambda)$, we would usually use the methods (\ref{PQR}) from \cite{FPP2} 
because they require knowledge of only the first few derivatives of $q(x)$.
But when $q$ has the required special forms, the method (\ref{fPQR})
from \cite{FPP4}, often more efficient, can also be used. 

\section{Numerical Estimation of the Spectral Density Function $f(\lambda)$}   %  6
\setcounter{equation}{0}

In this section we test our implementation of the various numerical methods
from the previous section on the examples listed in section 6, for which exact 
formulas are known for the spectral density function. Then we also test  a more interesting example from 
quantum chemistry for which exact formulas are not known.

For the hydrogen atom potential, Example~1 (equation (\ref{6.1})) with $a = 1$, 
Table~8.1 has numerical output when $\ell = 1$ for many $\lambda$ with the 
methods $F_x^1$, $F_x^2$, $F_x^3$, and $f_x^6$ using the notation of the 
previous section.  Table~8.2 shows the analogous data when $\ell = 2$.
Shown are only the errors: absolute when the answer is less than one and relative
otherwise.  A tolerance of $10^{-14}$ was used for the numerical integration 
of the initial value problem for (\ref{6.1}), starting at $x_0(\lambda)$ from (\ref{7.2}).
Consequently, table entries this
small represent errors in $y$ as well as errors due to finite $x$.  Note 
that for a fixed accuracy, generally larger matching points $x$ are needed 
when $\lambda$ is smaller.   As expected, the higher order methods $F^3$ 
and $f^6$ are much superior. 

\pagebreak
\vskip 0.2in
\begin{center}
Table 8.1. Numerical Error: H Atom ($\ell=1$) $q(x) = -1 / x + 2 /x^2$.
\vskip 4pt

\begin{tabular}{c|r@{}l r@{}l r@{}l r@{}l r@{}l}
$\lambda$&\multicolumn{2}{c}{$x=x(\lambda)$}&\multicolumn{2}{c}{$F^1_x$}&
\multicolumn{2}{c}{$F^2_x$}&\multicolumn{2}{c}{$F_x^3$}&
\multicolumn{2}{c}{$f^6_x$}\\
\hline
  0.1 & 320&.0 & $ $4&.60($-$4) & $-$2&.10($-$8) & $ $5&.62($-$12) & $ $2&.13($-$14) \\
  0.2 & 225&.0 & $ $2&.19($-$5) & $-$1&.18($-$9) & $ $4&.28($-$13) & $ $4&.48($-$14) \\
  0.4 & 160&.0 & $ $5&.08($-$4) & $-$2&.34($-$8) & $ $6&.36($-$12) & $ $4&.56($-$15) \\
    1 & 100&.0 & $ $5&.53($-$4) & $-$2&.62($-$8) & $ $7&.49($-$12) & $ $1&.25($-$13) \\
    2 &  71&.0 & $ $9&.52($-$4) & $-$4&.39($-$8) & $ $1&.21($-$11) & $ $9&.60($-$14) \\
    4 &  50&.0 & $-$2&.55($-$4) & $ $1&.22($-$8) & $-$2&.95($-$12) & $ $4&.89($-$13) \\
   10 &  32&.0 & $ $1&.41($-$3) & $-$5&.93($-$8) & $ $1&.48($-$11) & $ $3&.25($-$13) \\
   20 &  22&.5 & $-$2&.79($-$4) & $ $1&.14($-$8) & $-$1&.81($-$12) & $ $7&.83($-$13) \\
   40 &  16&.0 & $ $3&.52($-$4) & $-$1&.21($-$8) & $ $2&.91($-$12) & $ $8&.16($-$13) \\
  100 &  10&.0 & $-$3&.39($-$4) & $ $8&.44($-$9) & $ $2&.52($-$13) & $ $2&.61($-$13) \\
  200 &   7&.0 & $ $2&.32($-$4) & $-$2&.34($-$9) & $-$1&.91($-$12) & $ $3&.02($-$13) \\
  400 &   5&.0 & $-$1&.04($-$4) & $-$1&.66($-$9) & $ $2&.78($-$12) & $ $2&.60($-$13) \\
 1000 &   3&.2 & $-$4&.06($-$6) & $-$6&.53($-$10)& $ $7&.80($-$13) & $ $2&.82($-$13) \\
 2000 &   2&.2 & $ $5&.84($-$6) & $ $5&.12($-$9) & $-$2&.91($-$12) & $ $2&.86($-$13) \\
 4000 &   1&.6 & $ $3&.49($-$6) & $-$1&.79($-$9) & $ $1&.34($-$12) & $ $2&.84($-$13) \\
10000 &   1&.0 & $ $2&.67($-$5) & $-$6&.60($-$9) & $ $3&.81($-$12) & $ $2&.95($-$13) \\
\end{tabular}
\end{center}
\vskip 14pt

%\pagebreak
\vskip 0.2in
\begin{center}
Table 8.2.  Numerical Error: H Atom ($\ell=2$) $q(x) = -1 / x + 6 /x^2$.
\vskip 4pt

\begin{tabular}{c|r@{}l r@{}l r@{}l r@{}l r@{}l}
$\lambda$&\multicolumn{2}{c}{$x=x(\lambda)$}&\multicolumn{2}{c}{$F^1_x$}&
\multicolumn{2}{c}{$F^2_x$}&\multicolumn{2}{c}{$F_x^3$}&
\multicolumn{2}{c}{$f^6_x$}\\
\hline
  0.1 & 320&.0 & $-$3&.13($-$6) & $ $1&.35($-$10)& $-$3&.45($-$14) & $ $3&.95($-$15) \\
  0.2 & 225&.0 & $-$4&.79($-$6) & $ $2&.13($-$10)& $-$5&.60($-$14) & $ $2&.98($-$15) \\
  0.4 & 160&.0 & $-$4&.28($-$6) & $ $1&.83($-$10)& $-$4&.54($-$14) & $ $3&.00($-$15) \\
    1 & 100&.0 & $-$2&.82($-$5) & $ $1&.20($-$9) & $-$2&.98($-$13) & $ $8&.32($-$15) \\
    2 &  71&.0 & $-$7&.24($-$5) & $ $2&.87($-$9) & $-$6&.41($-$13) & $ $1&.92($-$14) \\
    4 &  50&.0 & $ $8&.09($-$5) & $-$2&.97($-$9) & $ $5&.83($-$13) & $ $5&.28($-$15) \\
   10 &  32&.0 & $-$8&.59($-$4) & $ $2&.24($-$8) & $-$7&.48($-$13) & $ $1&.69($-$13) \\
   20 &  22&.5 & $ $2&.80($-$4) & $-$3&.67($-$9) & $-$1&.75($-$12) & $ $8&.68($-$13) \\
   40 &  16&.0 & $-$2&.36($-$4) & $-$2&.06($-$9) & $ $4&.70($-$12) & $ $6&.52($-$14) \\
  100 &  10&.0 & $ $1&.69($-$4) & $ $1&.66($-$8) & $-$1&.22($-$11) & $ $8&.75($-$14) \\
  200 &   7&.0 & $-$4&.57($-$5) & $-$2&.62($-$8) & $ $1&.71($-$11) & $ $3&.02($-$14) \\
  400 &   5&.0 & $-$3&.42($-$5) & $ $2&.25($-$8) & $-$1&.26($-$11) & $-$5&.91($-$14) \\
 1000 &   3&.2 & $-$8&.80($-$6) & $ $2&.46($-$9) & $-$6&.75($-$13) & $ $8&.21($-$13) \\
 2000 &   2&.2 & $ $1&.06($-$4) & $-$2&.25($-$8) & $ $1&.24($-$11) & $ $8&.62($-$13) \\
 4000 &   1&.6 & $-$2&.78($-$5) & $ $5&.40($-$9) & $-$2&.15($-$12) & $ $8&.09($-$13) \\
10000 &   1&.0 & $-$1&.42($-$4) & $ $2&.38($-$8) & $-$1&.17($-$11) & $ $8&.94($-$13) \\
\end{tabular}
\end{center}
\vskip 14pt

The next choice is the Bessel Equation, Example~3 (equation (\ref{6.9})), with $\nu = 1/3$:
\begin{equation}   \label{BESS3}          %  760a
   q(x) = -5 / (36 x^2). 
\end{equation}
The error behavior is similar to that in the previous tables.

\pagebreak
\vskip 0.2in
\begin{center}
Table 8.3.  Numerical Error:  $q(x) = -5 / (36 x^2)$ 
\vskip 4pt
\begin{tabular}{c|r@{}l r@{}l r@{}l r@{}l r@{}l}
$\lambda$&\multicolumn{2}{c}{$x=x(\lambda)$}&\multicolumn{2}{c}{$F^1_x$}&
\multicolumn{2}{c}{$F^2_x$}&\multicolumn{2}{c}{$F_x^3$}&
\multicolumn{2}{c}{$f^6_x$}\\
\hline
  0.1 & 320&.0 & $-$3&.26($-$7) & $ $4&.61($-$11)& $ $3&.21($-$14) & $ $5&.37($-$14) \\
  0.2 & 225&.0 & $ $1&.22($-$6) & $-$1&.13($-$10)& $ $1&.59($-$13) & $ $6&.86($-$14) \\
  0.4 & 160&.0 & $-$5&.18($-$7) & $ $7&.31($-$11)& $ $4&.90($-$14) & $ $8&.34($-$14) \\
    1 & 100&.0 & $ $2&.34($-$6) & $-$3&.48($-$10)& $ $2&.85($-$13) & $ $1&.12($-$13) \\
    2 &  71&.0 & $ $3&.30($-$6) & $-$4&.93($-$10)& $ $3&.92($-$13) & $ $1&.47($-$13) \\
    4 &  50&.0 & $ $3&.71($-$6) & $-$5&.53($-$10)& $ $4&.38($-$13) & $ $1&.63($-$13) \\
   10 &  32&.0 & $-$1&.51($-$6) & $ $2&.14($-$10)& $ $1&.49($-$13) & $ $2&.49($-$13) \\
   20 &  22&.5 & $ $5&.26($-$6) & $-$7&.86($-$10)& $ $6&.76($-$13) & $ $2&.85($-$13) \\
   40 &  16&.0 & $-$1&.78($-$6) & $ $2&.51($-$10)& $ $1&.58($-$13) & $ $2&.76($-$13) \\
  100 &  10&.0 & $ $5&.92($-$6) & $-$8&.26($-$10)& $ $7&.22($-$13) & $ $2&.83($-$13) \\
  200 &   7&.0 & $-$5&.92($-$6) & $ $9&.12($-$10)& $-$1&.74($-$13) & $ $2&.95($-$13) \\
  400 &   5&.0 & $ $5&.92($-$6) & $-$8&.82($-$10)& $ $7&.36($-$13) & $ $2&.97($-$13) \\
 1000 &   3&.2 & $-$1&.78($-$6) & $ $2&.51($-$10)& $ $1&.54($-$13) & $ $2&.72($-$13) \\
 2000 &   2&.2 & $-$5&.78($-$6) & $ $8&.89($-$10)& $-$1&.65($-$13) & $ $2&.91($-$13) \\
 4000 &   1&.6 & $-$1&.78($-$6) & $ $2&.51($-$10)& $ $1&.70($-$13) & $ $2&.88($-$13) \\
10000 &   1&.0 & $ $5&.92($-$6) & $-$8&.82($-$10)& $ $5&.98($-$12) & $ $2&.63($-$13) \\
\end{tabular}
\end{center}
\vskip 14pt

In the papers of Bain et al \cite{BAIN}, Br\"andas et al \cite{BRA1}, and
Engdahl et al \cite{ENG1}, \cite{ENG2}, \cite{ENG3} can be found a potential 
with a ``barrier'' near $x = 2$ and decaying rapidly to zero as 
$x \to \infty$:
\begin{equation}     \label{BARRIER}
    q(x) = \sfrac{\ell (\ell+1)}{x^2} - \sfrac{a}{x} + 15 x^2 e^{-x}.
\end{equation}
For this example we have $q_{0} = \ell (\ell+1)$, $q_{1} = -a$,
$q_2 = q_3 = 0$, and for $k \ge 2$
\[  q_{k+2} = (-1)^k \sfrac{15}{(k-2)!}.   \]
For $a = 1$, Table~8.4 displays the results for 
several values of $x$ to show the rapid convergence as $x \to \infty$.
A tolerance of $10^{-14}$ was used for the numerical integration of $y$.

%\pagebreak
\vskip 0.2in
\begin{center}
Table 8.4. $F_x^3$ Estimates for Barrier Potential (\ref{BARRIER}) with $a = 1$.
\vskip 4pt

\begin{tabular}{c|ccc|ccc}
 $x$  & $\ell = 0$ & $\ell = 1$ &  $\ell = 2$ & $\ell = 0$ & $\ell = 1$ &  $\ell = 2$ \\
\hline
      &            &$\lambda = 7$&             &            &$\lambda=10$&              \\
   5. & 0.142809355& 0.019804657& 0.0004166628 & 1.686525464& 1.728228916& 0.1242771047 \\
  10. & 0.142828980& 0.019801387& 0.0004162081 & 1.686646374& 1.728086680& 0.1242724151 \\
  15. & 0.142829143& 0.019801395& 0.0004162075 & 1.686647533& 1.728085796& 0.1242722883 \\
  20. & 0.142829149& 0.019801396& 0.0004162075 & 1.686647559& 1.728085772& 0.1242723271 \\
  25. & 0.142829149& 0.019801396& 0.0004162075 & 1.686647559& 1.728085772& 0.1242723271 \\
      &            &$\lambda = 20$&            &            &$\lambda=40$&              \\
   5. & 1.999374819& 4.314112367& 2.8784054457 & 2.558971562& 11.31991454& 17.270757428 \\
  10. & 1.999374881& 4.314112204& 2.8784044730 & 2.558971266& 11.31991563& 17.270756351 \\
  15. & 1.999374882& 4.314112311& 2.8784043230 & 2.558971293& 11.31991552& 17.270756523 \\
  20. & 1.999374882& 4.314112307& 2.8784043239 & 2.558971293& 11.31991552& 17.270756528 \\
  25. & 1.999374882& 4.314112307& 2.8784043240 & 2.558971293& 11.31991552& 17.270756528 \\
\end{tabular}
\end{center}
\vskip 14pt

\section{Numerical Estimation of the Spectral Function $\rho(\lambda)$}  % 9
\setcounter{equation}{0}

Associated with the density function $f(\lambda)$ is the associated spectral function 
defined by
\begin{equation}    \label{RHODEF}
  \rho(\lambda): = \int_0^{\lambda} f(\mu) \, d \mu.  
\end{equation}
In \cite{FPP1} for problems regular at $x=0$ we estimated $\rho$ using 
$F^1$ and compared with the package SLEDGE \cite{SLEDGE,TOMS98}; generally the $\rho(\lambda)$
computation using a quadrature routine for (\ref{RHODEF}) and the $F^1$ formula ran 
considerably faster than SLEDGE, but still had the drawback that rather large $x-$intervals
were required for the $F^1$ calculation.  Here we apply 
the methods of this paper for computing $\rho$ by  estimating $f$ in (\ref{RHODEF}), and performing a
quadrature, to the examples in section 6 for which exact answers are known, and again compare with SLEDGE.

The SLEDGE software for estimating $\rho(\lambda)$ is based on the 
Levitan-Levinson characterization of the spectral function as a limit of step spectral functions 
over a finite interval approximation (second formula in (\ref{5.8})); this is a 
totally different approach than the present approach of this paper which relies on the 
family of $F^{j}$-approximants, together with the quadrature in (\ref{RHODEF}).  
%Its much greater computing times for a given accuracy were also observed 
%in \cite{FPP1} for problems regular at $x=0$.
For the case of two singular endpoints, the performance of SLEDGE for computing the
spectral function on examples having explicit formulas for the spectral function
was reported on in \cite{TOMS98}. As reported there, one of the major weaknesses 
of the SLEDGE package is obtaining high accuracy in the $\rho(\lambda)$ calculation
when $\lambda$ is large;  this is due primarily to the fact that SLEDGE does not 
rely on asymptotic approximations for the eigenvalues and eigenfunction norm reciprocals,
but computes them numerically as required for implementing the $\rho_b(\lambda)$-formula in (\ref{5.8}).
Experience in using SLEDGE on doubly singular problems is that very large computing
times are required due to the computation of large numbers of eigenvalue - eigenfunction norm
pairs, and that there is significant loss of accuracy when $\lambda$ becomes sufficiently large.
As the timing and accuracy data of this section shows, doubly singular problems can be handled
with high accuracy and much reduced computing times by making use of the $F_x^j$-approximants and the $f_x^N$-approximants
of this paper, along with the quadratures for computing $\rho(\lambda)$ in  (\ref{RHODEF}); this represents a major
improvement in computational technique over the SLEDGE algorithm for spectral function computation.

%  The computation time of
%SLEDGE is greatly reduced at modest requested tolerances when SLEDGE can 
%use its built-in asymptotic formulas for large $\lambda$.   While these 
%are implemented for many cases of Sturm-Liouville problems, they do not 
%include the doubly singular cases of this paper.  Consequently, it should 
%be no surprise that SLEDGE suffers in comparison with the methods of 
%Section~4 when $\lambda$ is large.  

Following SLEDGE, we assume approximations are sought for a finite  set of $\lambda$-values
in the continuous spectrum,  $(0,\infty)$,
ordered so that
\[  0  < \lambda_1 < \lambda_2 < ... < \lambda_m.  \]
Then with $\rho(0)$ given (or computed via SLEDGE) we estimate
\[  \rho(\lambda_1) = \rho(0) + \int_0^{\lambda_1} f(\mu) \, d \mu \]
and
for $j = 2, 3, \ldots, m$
\begin{equation}    \label{RHOSUM}
  \rho(\lambda_j) = \rho(\lambda_{j-1}) + 
    \int_{\lambda_{j-1}}^{\lambda_j} f(\mu) \, d \mu   
\end{equation}
using an adaptive quadrature code with $f$ replaced by $F_x^j$ or $f_x^N$ approximations.
Here we report some timing and accuracy data for the four examples listed in section 6.

%\begin{enumerate}
%\item {\bf{ Example 9, Bessel Equa with $\nu=1$:}} \hspace{0.1in} $ q(x) = 0.75 / x^2 $      
%\item  {\bf{Example 1, H-atom with $\ell = 1$:}} \hspace{0.1in} $ q(x) = -1 / x + 2 / x^2 $
%\item  {\bf{Example 2, Repulsive Coulomb with $\ell = 1$:}} \hspace{0.1in} $ q(x) = 1 / x + 2 / x^2 $
%\item {\bf{ Example 9, Bessel Equa with $\nu=\sfrac{1}{3}$:}} \hspace{0.1in}  $q(x) = -5 / (36 x^2) $      
%\end{enumerate}

The spectral functions on $(0,\infty)$ for these examples are known in closed form by putting the exact spectral
density functions from (\ref{6.4}), (\ref{6.8}), (\ref{6.12}), and (\ref{6.16}) into (\ref{RHODEF}) and performing
an exact integration. The resulting closed form formulas for $\rho(\lambda)$ were used for the four  examples to compare with the
numerical approximations and to generate the `exact' error; the errors are taken as absolute if the exact value
of $\rho(\lambda)$ is less than one, and relative otherwise.
       
%The first example is one with known closed form solution -- this will
%provide a check for the algorithm and code.   The potential is from the
%Bessel equation of order one: 
%\begin{equation}   \label{BESS1}
%   q(x) = 0.75 / x^2. 
%\end{equation}
%On the interval $(0, \infty)$ this is singular at each endpoint: 
%at $x =0 $ it is limit point and nonoscillatory.  It is known that 
%$\rho(\lambda) = \lambda^2 / 16$. The potential is of the form (\ref{GENq}), 
%and we have chosen $N=7$. 

For the Bessel equation of order 1, Example 4 (equation (\ref{6.13} with N=1),  we used $N=7$ in the
approximation (\ref{fPQR})  (that is, the scheme from \cite[Sec 4]{FPP4}).  Output data for six $\lambda$-values is
displayed in  Table~9.1.   The quadrature tolerance was $10^{-8}$ and the tolerance for the initial value problem was $10^{-9}$. 
Note that, at this tolerance, all the apparent error arises from the first 
integration interval and is passed on through the sum in (\ref{RHOSUM}).  

\vskip 0.2in
\begin{center}
Table 9.1. {\bf $\rho(\lambda)$ for the first order Bessel equation on $(0, \infty)$} {\bf(q(x) = $0.75 / x^2$)}   % 740a
\vskip 4pt
\begin{tabular}{c|r@{}l r@{}l r@{}l  r@{}l  r@{}l  r@{}l}
$x$&\multicolumn{2}{c}{$\lambda = 1$}&\multicolumn{2}{c}{$\lambda = 2$}
&\multicolumn{2}{c}{$\lambda = 4$}&\multicolumn{2}{c}{$\lambda = 10$}
&\multicolumn{2}{c}{$\lambda = 20$}&\multicolumn{2}{c}{$\lambda = 40$}\\
\hline
 6. &0&.06318042&0&.25068042&1&.00068041&6&.25068041&25&.00068041&100&.00068041\\
12. &0&.06254252&0&.25004252&1&.00004252&6&.25005252&25&.00004252&100&.00004252\\
24. &0&.06250266&0&.25000266&1&.00000266&6&.25000266&25&.00000266&100&.00000266\\
36. &0&.06250000&0&.25000000&1&.00000000&6&.25000000&25&.00000000&100&.00000000\\
\end{tabular}
\end{center}
\vskip 14pt

The data for estimating $f(\lambda)$ in section 8 showed that as $\lambda$ 
gets larger, smaller values of the matching point $x=x(\lambda)$ are needed for a given accuracy.  The
data in Table~9.1 for $\rho$  exhibit this  phenomena.  Hence,
from efficiency considerations in order to compute $\rho(\lambda)$ we want 
the choice of matching point $x$ to vary with the integer $N$, or even better, pointwise 
with $\lambda$.   For $q(x) = 0.75 / x^2$ it can be shown that the absolute
error in $f(\lambda)$ for a given $N$ is proportional to   
\[    \sfrac{1}{x^{N+2} \lambda^{N-1/2}}.  \]
This suggests that taking $x \sim 1 / \lambda^{(2N-1)/(2N+4)}$ would be a 
good heuristic for the matching point.  Similarly, for the family of approximants, $F^N$, from (\ref{PQR})--(\ref{PQR3}) the 
corresponding form for the absolute error in $f(\lambda)$ is
\[    \sfrac{1}{x^{2N} \lambda^{N-1/2}},  \]
so that $x \sim 1 / \lambda^{(2N-1) / (4N)}$ would be an appropriate matching point.  The
latter is roughly $1 / \sqrt{\lambda}$.  As mentioned earlier, these 
formulas would change for a different $q$.  Also, the heuristic would
differ for relative errors.  A similar analysis for the general potential
$q(x) = A /x + B / x^2$ suggests a good first value of matching point would be  
\begin{equation}   \label{xHEUR}
  x = x(\lambda) := |A| / (2 \lambda) + \sqrt{A^2 / (4 \lambda^2) + |B| / \lambda}.  
\end{equation}

We have written a research code, called AutoB, for which the only inputs 
required are the set of $\lambda$ points, the $\rho(0)$ value, and the 
accuracy desired.   If $q$ has the form of (\ref{GENq}) or the form of similar potentials in
\cite{FPP4}, then we use the appropriate $f_{x}^N$ formula in (\ref{fPQR})
with $N$ chosen to be a function of the accuracy sought.  Otherwise, we 
use $F^3$ from (\ref{PQR}) which requires knowledge of derivatives of $q$.
We report the performance of AutoB on the four examples in section 6
using (\ref{xHEUR}) as the initial choice of matching point.   Given a 
prescribed tolerance $\tau$, for each $\lambda$ the matching point $x$ is then
increased until
\[ |\mbox{error estimate}| \le \max\{1, |\mbox{output value}|\} \; \tau  \]
holds.
%AutoB has been tested on several examples, results from four of which we
%summarize in Table~7.2.  The first line is for potential (\ref{BESS1}),
%the second is the hydrogen atom ($\ell = 1$) potential
%\begin{equation}   \label{HATMx}          %  710b
%   q(x) = -1 / x + 2 / x^2,
%\end{equation}
%the third is of the form of Example~5 in Sec.~3 
%\begin{equation}   \label{BESS2}          %  735b
%   q(x) = 1 / x + 2 / x^2, 
%\end{equation}
%and the final potential is (\ref{BESS3}) from the previous section. 
Estimates at various $\tau$ were sought for the following set of sixteen $\lambda$ values:
\begin{equation}  \label{RhoSet}
  \{ 0.1, 0.2, 0.4, 1, 2, 4, 10, 20, 40, 100, 200, 400, 1000, 2000, 4000, 10000 \}.   
\end{equation}
The error shown is the maximum (relative when the `exact' $\rho > 1$, absolute otherwise)
over the set of sixteen $\lambda$ values.
For many of these runs the heuristics were overly conservative, but the times 
are nevertheless quite small.  For the Bessel equation of order $\sfrac{1}{3}$, 
Example 3 (equation (\ref{6.9}) with $\nu=1/3$), the relatively large computing times 
were due to difficulties near $\lambda=0$.

\vskip 0.2in
\begin{center}
Table 9.2.  {\bf{AutoB results for several tolerances and Four Potentials on (0,$\infty$).}}
\vskip 4pt
\begin{tabular}{lcccccccc}
&\multicolumn{2}{c}{$\tau=10^{-4}$}&\multicolumn{2}{c}{$\tau=10^{-6}$}&
\multicolumn{2}{c}{$\tau=10^{-8}$}&\multicolumn{2}{c}{$\tau=10^{-10}$}  \\
Potential&error&time &   error  & time &   error  & time &   error  & time    \\
\hline
1. Ex4($\nu$=1) &1.30($-$7)& 0.28 &1.50($-$9)& 0.50 &1.92($-$12)& 1.06 &1.73($-$13)&  2.46 \\  % 740a
2. Ex1($\ell$=1) &1.78($-$6)& 0.42 &7.88($-$7)& 0.84 &1.50($-$8) & 3.11 &7.93($-$11)& 18.50 \\  % 710b
3. Ex2($\ell$=1) &2.49($-$6)& 0.37 &3.37($-$7)& 0.74 &8.54($-$9) & 3.95 &1.10($-$10)& 20.97 \\  % 735b
4. Ex3($\nu$=1/3) &9.57($-$5)& 0.49 &1.04($-$6)& 3.00 &1.52($-$8) &41.95 &1.13($-$10)&894.77     % 760a
\end{tabular}
\end{center}
\vskip 14pt

For comparison, output from the SLEDGE program  is shown for the Bessel equation of order Example 4 with $\nu=1$, Example 4 (equation (\ref{6.13}) in Table~9.3.  Since SLEDGE is known \cite{TOMS98} to have difficulty with large values of 
$\lambda$, we ran the program at various choices of $\tau$ only on the first $n$ $\lambda$-values
from (\ref{RhoSet}) with $n = 1, 2, \ldots, 9$.  The final line ($n = 16$)
is output from AutoB using all sixteen $\lambda$ values up to $\lambda=10000$.  Clearly, AutoB is 
much more reliable than SLEDGE. Similar results were observed for other doubly singular potentials.

%\vskip 0.2in
\begin{center}
Table 9.3. {\bf SLEDGE output for the Bessel equation of order 1}{\bf (q(x) = $0.75 / x^2$)}   % 740a
\vskip 4pt
\begin{tabular}{ccccccccc}
&\multicolumn{2}{c}{$\tau=10^{-3}$}&\multicolumn{2}{c}{$\tau=10^{-4}$}&
\multicolumn{2}{c}{$\tau=10^{-5}$}&\multicolumn{2}{c}{$\tau=10^{-6}$}  \\
$n$&   error  & time &   error  & time &   error  & time &   error  & time \\
\hline
 1 &1.31($-$4)& 0.02 &3.87($-$5)& 0.09 &9.74($-$6)& 0.50 &9.92($-$6)& 1.58 \\
 2 &6.86($-$4)& 0.02 &1.24($-$5)& 0.23 &1.07($-$5)& 1.31 &9.87($-$6)& 8.43 \\
 3 &2.72($-$4)& 0.14 &1.32($-$5)& 0.34 &2.95($-$5)& 2.34 &6.70($-$6)&$>$96.53\\
 4 &2.71($-$4)& 0.16 &1.79($-$4)& 0.50 &3.20($-$5)& 8.45 &6.70($-$6)&$>$177.28\\
 5 &7.16($-$4)& 0.19 &1.79($-$4)& 0.70 &3.20($-$5)&16.40 &          &      \\
 6 &1.92($-$3)& 0.21 &1.79($-$4)& 1.71 &3.20($-$5)&29.51 &          &      \\
 7 &7.04($-$3)& 0.96 &1.83($-$4)& 7.67 &8.65($-$5)&51.64 &          &      \\
 8 &1.47($-$2)& 1.83 &1.83($-$4)&11.39 &2.18($-$4)&75.87 &          &      \\
 9 &2.23($-$2)&$>$12.00&3.31($-$4)&20.48&2.83($-$4)&$>$465.99&      &      \\
\hline
AutoB \\
16 &3.69($-$6)& 0.18 &1.30($-$7)& 0.28 &5.74($-$9)&  0.36&1.50($-$9)& 0.50 \\
\end{tabular}
\end{center}
\vskip 14pt

To illustrate the superiority of the new code AutoB over SLEDGE we also ran comparisons
on timing and accuracy the Hydrogen Atom potential with $\ell=1$, Example 1 (equation (\ref{6.1}).  The output values
for $\rho(\lambda)$ obtained for each of the sixteen $\lambda$ values in (\ref{RhoSet})
for four choices of the tolerance levels are displayed in Table~9.4. The $>$ in the time
needed for SLEDGE indicates that it stopped (too much time) before the user requested input accuracy was achieved.
Since SLEDGE could not achieve $10^{-3}$-accuracy over the whole range of $\lambda$-values, 
no SLEDGE runs for tighter tolerances are listed. 
%As an illustration of the performance of the new codes produced under this grant, the following table from \cite{REF8} compares
%SLEDGE with the new codes on time and accuracy for computing the spectral function (\ref{EQ3_14}) associated with the Hydrogen
%problem (\ref{EQ3_26})-(\ref{EQ3_27}) for $\ell=1$; here the above formula (\ref{EQ3_29}) is used to generate the `exact' answers for the
%accuracy comparisons. 
As the data shows, SLEDGE has much difficulty to compute highly accurate results for large values of 
$\lambda$, while the new codes are capable of quite high accuracy in much less computing time. Similar testing for the
Hydrogen Atom potential using the $F_x^j$ approximants was also done in the thesis of Mark Schuster \cite{MARK}.

\vskip 0.2in
\begin{center}
Table 9.4: {\bf{Comparison of SLEDGE with AutoB for Hydrogen problem with $\ell = 1.$}}
\vskip 4pt
\begin{tabular}{cc|ccccc}
$\lambda$& Exact &  SLEDGE  &  AutoB  &  AutoB  &  AutoB   &   AutoB  \\
\hline
 0.1&0.005621362  &0.0056   &0.0056   &0.00562  &0.0056214 & 0.005621362 \\
 0.2&0.010067470  &0.0100   &0.0100   &0.01007  &0.0100675 & 0.010067470 \\
 0.4&0.022334469  &0.0222   &0.0223   &0.02233  &0.0223345 & 0.022334470 \\
 1.0&0.087358065  &0.0867   &0.0874   &0.08736  &0.0873581 & 0.087358074 \\
 2.0&0.298717032  &0.2966   &0.2987   &0.29872  &0.2987170 & 0.298717047 \\
 4.0&1.166166722  &1.1577   &1.1662   &1.16617  &1.1661667 & 1.166166736 \\
 10.&8.206942681  &8.1493   &8.2069   &8.20694  &8.2069402 & 8.206942643 \\
 20.&38.98117554  &38.691   &38.981   &38.9812  &38.981169 & 38.98117542 \\
 40.&194.5884791  &192.80   &194.59   &194.588  &194.58845 & 194.5884785 \\
100.&1719.215348  &1706.1   &1719.2   &1719.21  &1719.2140 & 1719.215343 \\
200.&9188.295022  &9068.2   &9188.3   &9188.29  &9188.2922 & 9188.294986 \\
400.&49923.13741  &49137.   &49923.   &49923.1  &49923.126 & 49923.13720 \\
1000.&475962.2250 &462178.  &475962.  &475961.  &475962.23 & 475962.2484 \\
2000.&2644112.132 &2479326. &2644120. &2644111. &2644111.6 & 2644112.120 \\
4000.&14766762.30 &13806576.&14766758.&14766759.&14766759. & 14766762.23 \\
10000.&144274264.9&112742503&144274122&144274122&144274123 & 144274264.3 \\
\hline
time (sec)&       &$>$177.5 &  0.31   &  0.42   &   0.83   &   3.75    \\ 
 RelErr &         &$10^{-3}$&$10^{-3}$&$10^{-4}$& $10^{-6}$& $10^{-8}$ \\
 AbsErr &         &$10^{-3}$&$10^{-3}$&$10^{-4}$& $10^{-6}$& $10^{-8}$ \\
\end{tabular}
\end{center}
\vskip 14pt

{\bf Remark.} High accuracy in the spectral function computation for the Bessel
equation on $(0,\infty)$, Examples 3 and 4, was also achieved in \cite[Sec 4]{TOMS98}; this, however,
was done by inserting asymptotic formulas for the eigenvalues and eigenfunction norm reciprocals
for the Bessel equation on $(0,b]$ into the SLEDGE code (bypassing the SLEDGE computation of these
quantities); but, of course, this was not an automatic procedure applicable to other problems with two singular endpoints.

\end{document}